\documentclass[12pt,a4paper]{amsart}

\usepackage{fullpage}
\usepackage[pdftex]{epsfig}

\usepackage{amssymb}  
\usepackage{latexsym} 
\usepackage{comment}
\usepackage{color}
\usepackage{enumerate}
\usepackage{amsthm}
\usepackage[all]{xy}

\usepackage{tikz}
\usetikzlibrary{matrix,arrows,automata,snakes}
\usepackage[textwidth=20mm,textsize=small]{todonotes}
\usepackage{charter,euler} 
\DeclareSymbolFont{bchoperators}{T1}{bch}{m}{n}
\SetSymbolFont{bchoperators}{bold}{T1}{bch}{b}{n}
\makeatletter
\renewcommand{\operator@font}{\mathgroup\symbchoperators}
\makeatother

\usepackage{titlesec} 
\titleformat{\section}{\normalfont\bfseries\filcenter}{\thesection}{1em}{}
\titleformat{\subsection}{\normalfont\bfseries}{\thesubsection}{1em}{}
\titleformat{\subsubsection}{\normalfont\bfseries}{\thesubsubsection}{1em}{}

\numberwithin{equation}{section}

\definecolor{darkgreen}{rgb}{0,0.5,0}
\definecolor{rem}{rgb}{0.8,0,0}
\definecolor{new}{rgb}{0.7,0,0.6}
\definecolor{reply}{rgb}{0,0,0.8}
\usepackage[
        colorlinks, citecolor=darkgreen,
        backref,
        pdfauthor={Maarten Derickx, Sheldon Kamienny, William Stein, Michael Stoll},
        pdftitle={Torsion points on elliptic curves over number fields of small degree},
]{hyperref}
\usepackage[capitalize]{cleveref}
\usepackage[alphabetic,backrefs,lite]{amsrefs} 

\newtheorem{theorem}{Theorem}[section]
\newtheorem{lemma}[theorem]{Lemma}
\newtheorem{proposition}[theorem]{Proposition}

\newtheorem{corollary}[theorem]{Corollary}

\theoremstyle{definition}
\newtheorem{definition}[theorem]{Definition}

\theoremstyle{remark}
\newtheorem*{remark}{Remark}

\renewcommand{\O}{\mathcal{O}} 
\newcommand{\Ohat}{\hat{\text{\rule{0mm}{1.75ex}}\O}} 

\newcommand{\set}[1]{\left\lbrace #1 \right\rbrace}
\newcommand{\diamondop}[1]{\langle #1 \rangle} 
\newcommand{\field}[1]{\mathbb{#1}}  
\newcommand{\Q}{\field{Q}} 

\newcommand{\R}{\field{R}} 
\newcommand{\C}{\field{C}} 
\newcommand{\Z}{\field{Z}} 
\newcommand{\F}{\field{F}} 
\newcommand{\T}{\field{T}} 
\renewcommand{\P}{\field{P}}

\newcommand{\calA}{\mathcal{A}}

\newcommand{\bfe}{\text{\bf e}}
\renewcommand{\frm}{\mathfrak{m}}
\newcommand{\HH}{\mathfrak{H}} 
\newcommand{\smm}[4]{\left(\begin{smallmatrix} #1 & #2 \\ #3 & #4 \end{smallmatrix}\right)}
\newcommand{\eps}{\varepsilon}

\newcommand{\ceil}[1]{\left\lceil #1 \right\rceil}
\newcommand{\floor}[1]{\left\lfloor #1 \right\rfloor}
\newcommand{\Bfloor}[1]{\Bigl\lfloor #1 \Bigr\rfloor}

\DeclareMathOperator{\Ann}{Ann}
\DeclareMathOperator{\CM}{CM}

\DeclareMathOperator{\cusps}{cusps}

\DeclareMathOperator{\End}{End}
\DeclareMathOperator{\Frob}{Frob}
\DeclareMathOperator{\Gal}{Gal}

\DeclareMathOperator{\gon}{gon}

\DeclareMathOperator{\im}{im}

\DeclareMathOperator{\ord}{ord}

\DeclareMathOperator{\primes}{Primes}
\DeclareMathOperator{\PSL}{PSL}
\DeclareMathOperator{\red}{red}
\DeclareMathOperator{\SL}{SL}

\DeclareMathOperator{\Tan}{Tan}
\DeclareMathOperator{\tors}{tors}
\DeclareMathOperator{\Ver}{Ver}

\begin{document}

\title[Torsion points on elliptic curves]%
      {Torsion points on elliptic curves \\ over number fields of small degree}

\author{Maarten Derickx}
\address{Mathematisch Instituut,
         Universiteit Leiden,
         P.O. Box 9512,
         2300 RA Leiden,
         The Netherlands}
\email{maarten@mderickx.nl}

\author{Sheldon Kamienny}
\address{University of Southern California,
         3620 South Vermont Ave., KAP 108,
         Los Angeles, California 90089-2532, USA}
\email{kamienny@usc.edu}

\author{William Stein}
\address{SageMath, Inc.,
         17725 SE 123rd Pl,
         Renton, WA 98059, USA.}
\email{wstein@gmail.com}

\author{Michael Stoll}
\address{Mathematisches Institut,
         Universit\"at Bayreuth,
         95440 Bayreuth, Germany.}
\email{Michael.Stoll@uni-bayreuth.de}

\date{January 14, 2021}

\keywords{Elliptic curve, torsion point, torsion subgroup, number fields of small degree}
\subjclass[2010]{Primary 11G05; Secondary 14G05, 14G25, 14H52}

\begin{abstract}
  We determine the set $S(d)$ of possible prime orders of $K$-rational points
  on elliptic curves over number fields~$K$ of degree~$d$,
  for $d = 4$, $5$, $6$, and~$7$.
\end{abstract}

\maketitle

\SelectTips{eu}{} 


\section{Introduction}

Let $K$ be an algebraic number field and let $E$ be an elliptic curve
over~$K$. Then the group~$E(K)$ of $K$-rational points on~$E$
is a finitely generated abelian group;
in particular, its torsion subgroup~$E(K)_{\tors}$ is a finite abelian
group, and one can ask which finite abelian groups can occur as the
torsion subgroup of~$E(K)$ for some elliptic curve over some number
field~$K$ of degree~$d$.

For $K = \Q$ (equivalently, $d = 1$), Mazur~\cites{mazur1,mazur2}
famously proved that the known finite list of possibilities
for the torsion subgroup is complete. This
was later extended by Merel~\cite{merel}, who
showed that for any given degree~$d$, there are only finitely
many possibilities for~$E(K)_{\tors}$ when $[K : \Q] = d$.

One key step in these finiteness results is to show that there
are only finitely many prime numbers~$p$ that can divide the order
of~$E(K)_{\tors}$, i.e., can occur as the order of an element
of~$E(K)$, for $K$ of degree~$d$. We therefore make the following
definition (following~\cite{kamiennymazur}).

\begin{definition}
  Let $d \ge 1$ be an integer. Then we define $S(d)$ to be the set
  of all prime numbers~$p$ such that there exists a number field~$K$
  of degree~$d$, an elliptic curve~$E$ over~$K$ and a point $P \in E(K)$
  such that $P$ has order~$p$.

  We write $\primes(x)$ for the set of all prime numbers~$p$
  such that $p \leq x$.
\end{definition}

Mazur showed that
\[ S(1) = \primes(7) . \]
Kamienny~\cite{kamienny2} determined
\[ S(2) = \primes(13) . \]
Merel~\cite{merel}*{Prop.~2 and~3} showed that
\[ S(d) \subseteq \primes(2^{d+1} d!^{5d/2}) \]
for $d \ge 4$. Parent~\cite{parent1} gave the better bound
(for all~$d$)
\[ S(d) \subseteq \primes(65 (3^d - 1) (2d)^6) . \]
However, Oesterl\'e had improved this already (as mentioned in Parent's
paper) to
\begin{equation} \label{eq:oesterle}
  S(d) \subseteq \primes((3^{d/2} + 1)^2)
\end{equation}
(except for not ruling out that $43 \in S(3)$)
in his unpublished notes~\cite{oesterle}.
It should be noted that Parent actually shows that his bound is valid
for prime \emph{power} order~$p^n$ of a torsion point when $p \ge 5$
(and he has similar bounds for powers of $2$ and~$3$); this is the main
point of his work.
Parent~\cites{parent2,parent3}, extending the techniques used by Mazur
and Kamienny and relying on Oesterl\'e's work, proved that
\[ S(3) = \primes(13) . \]

The main result of this paper is the following theorem,
which extends these results to $d = 4$, $5$, $6$, and~$7$.

\begin{theorem} \label{thm:main}
  \begin{align*}
    S(4) &= \primes(17), \\
    S(5) &= \primes(19), \\
    S(6) &= \primes(19) \cup \set{37}, \quad \text{and} \\
    S(7) &= \primes(23).
  \end{align*}
\end{theorem}

We also give a simplified proof of Parent's result on~$S(3)$.
Since we rely on Oesterl\'e's bound~\eqref{eq:oesterle}, a proof of which
has not been published so far, we include a proof here that is based
on Oesterl\'e's notes, which he kindly made available to us.

It is much easier to determine the set~$S'(d)$ of primes~$p$ such that
there are \emph{infinitely many} elliptic curves~$E$ over number fields~$K$
of degree~$d$ with distinct $j$-invariants that have a $K$-point of
order~$p$. This is mostly a question about the gonality of the modular
curve~$X_1(p)$. The following is known.

\begin{proposition}
  \begin{align*}
    S'(1) &= \primes(7), &
    S'(2) &= \primes(13), &
    S'(3) &= \primes(13), &
    S'(4) &= \primes(17), \\
    S'(5) &= \primes(19), &
    S'(6) &= \primes(19), &
    S'(7) &= \primes(23), &
    S'(8) &= \primes(23).
  \end{align*}
\end{proposition}

For $d = 1, 2, 3, 4$, this is shown in~\cites{mazur1,kamienny2,JKL1,JKL2},
respectively; for $5 \le d \le 8$, this follows
from~\cite{derickx_hoeij}*{Thm.~3}. Since clearly $S'(d) \subseteq S(d)$,
these results, together with the fact that a quadratic twist~$E_{6,37}$ over
the sextic number field \hbox{$K = \Q(\sqrt{5}, \cos(2\pi/7))$}
of the elliptic curve
\[ \text{\href{http://www.lmfdb.org/EllipticCurve/Q/1225/b/2}{$1225.b2$}}
    \colon y^2 + xy + y = x^3 + x^2 - 8x + 6
\]
has a point of order~$37$ over~$K$ \cite{elkies}*{Eq.~108}, reduce
the proof of~\cref{thm:main} to showing the inclusions ``$\subseteq$''.

We give the following more precise result in the case $d = 6$.

\begin{proposition} \label{prop:6_37}
  Let $K$ be a number field of degree~$6$ and let $E/K$ be an elliptic
  curve such that there is a point $P \in E(K)$ of exact order~$37$.
  Then $j(E) = j(E_{6,37}) = -9317$.
\end{proposition}

We prove~\cref{prop:6_37} at the end of~\cref{sec:ass_b}.

The gonality of~$X_1(p)$ grows like~$p^2$~\cite{abramovich}; this
implies that $S'(d) \subset \primes\bigl(O(\sqrt{d})\bigr)$. On the other
hand, denoting by $S_{\CM}(d)$ the set of primes that can occur as
orders of points on elliptic curves over a number field of degree~$d$
that have complex multiplication, the results of~\cite{CCS} show
that $S_{\CM}(s) \subset \primes\bigl(O(d)\bigr)$ and that
$3d+1 \in S_{\CM}(d)$ when $3d+1$ is prime. (Let $p = 3d+1$. There is a
pair of quadratic points defined over~$\Q(\sqrt{-3})$ with $j$-invariant
zero on~$X_0(p)$. The set-theoretic preimage gives a Galois orbit of points
of degree $2 \cdot \tfrac{p-1}{2} \cdot \tfrac{1}{3} = d$ on~$X_1(p)$,
since the covering $X_1(p) \to X_0(p)$ ramifies with index~$3$ above the
points with $j$-invariant zero.) So we will certainly have
$S'(d) \subsetneq S(d)$ for infinitely many~$d$. It is perhaps tempting
to assume that for large enough~$d$, the only sporadic points of degree~$d$
on~$X_1(p)$ are CM~points, as this seems to be the expectation for rational
points on modular curves. This would imply that $S(d) \subseteq \primes(3d+1)$
for large~$d$. However, consulting the table in~\cite{hoeij}, it appears
that there are many sporadic non-CM points (like the degree~$6$ points
on~$X_1(37)$ we have mentioned above). Still, the bound $p \le 3d+1$
is consistent with this information for $d \ge 13$.


\subsection*{The strategy}

To show the inclusions ``$\subseteq$'' in~\cref{thm:main}, we have
to verify that $p \notin S(d)$ for every prime number~$p$ that is not
in the set on the right hand side. This is equivalent
to the statement that all points of degree dividing~$d$ on the modular
curve~$X_1(p)$ over~$\Q$ are cusps. Recall that non-cuspidal points
on~$X_1(N)$, for $N \in \Z_{\ge 2}$, correspond to pairs~$(E,P)$, where
$E$ is an elliptic curve and $P \in E$ is a point of exact order~$N$.
See~\cref{sec:mod_curves} for some background on modular curves.

Now if $x \in X_1(p)(K)$ is a point
defined over a number field~$K$ of degree~$d$, but not over a smaller field,
then the sum of its Galois conjugates gives a $\Q$-rational effective
divisor of degree~$d$ on~$X_1(p)$. If $x$ is defined over a smaller
field~$K'$, then the degree~$d'$ of~$K'$ divides~$d$, and we can take
$d/d'$~times the sum of the conjugates of~$x$ to obtain a \hbox{$\Q$-rational}
effective divisor of degree~$d$ again. Effective divisors of degree~$d$
on a curve~$X$ correspond to points on its $d$th symmetric power~$X^{(d)}$
(which is the quotient of~$X^d$ by the natural action of the symmetric
group on $d$ letters). This leads to the following criterion.
We write $C_1(p)$ for the set of cusps on~$X_1(p)$.

\begin{lemma} \label{lem:criterion_alpha}
  Let $d \in \Z_{\ge 1}$ and let $p$ be a prime number.
  If the composition
  \[ \alpha \colon C_1(p)(\Q)^d \to X_1(p)(\Q)^d \to X_1(p)^{(d)}(\Q) \]
  of natural maps is surjective, then $p \notin S(d)$.

  If $p > 2d + 1$ and $p \notin S(d')$ for all $d' \le d$,
  then the map above is surjective.
\end{lemma}

\begin{proof}
  The assumption is equivalent to the statement that every $\Q$-rational
  effective divisor of degree~$d$ on~$X_1(p)$ is a sum of rational cusps.
  However, if there were a number field~$K$ of degree~$d$,
  an elliptic curve~$E$ over~$K$ and a point $P \in E(K)$ of order~$p$,
  then $(E, P)$ would give a $K$-rational non-cuspidal point on~$X_1(p)$
  and hence, by the discussion above, a $\Q$-rational effective divisor
  of degree~$d$ that is not supported on (rational) cusps, contradicting
  the assumption.

  For the converse, assume that the map is not surjective.
  Then there is a $\Q$-rational effective divisor~$D$ of degree~$d$ that is
  not supported on rational cusps. Since the irrational cusps on~$X_1(p)$
  form one Galois orbit of size $(p-1)/2 > d$, $D$ is not supported on
  cusps. This implies that there is a non-cuspidal point on~$X_1(p)$
  of degree $d' \le d$, hence $p \in S(d')$.
\end{proof}

We will follow the strategy that has been established in earlier work
by Mazur~\cite{mazur2}, Kamienny~\cites{kamienny1,kamienny2},
Merel~\cite{merel}, Oesterl\'e~\cite{oesterle} and
Parent~\cites{parent1,parent2,parent3}. We give an overview of the main
steps below; for a nice and more detailed account of Merel's proof
of the boundedness of~$S(d)$ for all~$d$ see~\cite{rebolledo}.

In our exposition, we refer to the existing literature for proofs
of many results we are using. Fairly detailed proofs of these statements
can be found in an earlier version of this paper~\cite{DKSSarXiv}
or in the doctoral thesis~\cite{DerickxThesis} of the first author.

The task is to show that $p \notin S(d)$ for $3 \le d \le 7$ and all
primes~$p$ not contained in the set on the right hand side of the equality
in~\cref{thm:main}. We use the criterion
of~\cref{lem:criterion_alpha}, in the equivalent form given below.
Before we formulate it, we make some definitions.

\begin{definition} \label{def:residue_class}
  Let $\ell$ be a prime. We write $\Z_{(\ell)}$ for the localization of~$\Z$
  at the prime ideal $(\ell) = \ell\Z$.

  Let $X$ be a scheme over~$\Z_{(\ell)}$. We denote the natural map
  $X(\Z_{(\ell)}) \to X(\F_\ell)$ by~$\red_\ell$.
  Let $\bar{x} \in X(\F_\ell)$. Then $\red_\ell^{-1}(\bar{x})$ is
  the \emph{residue class of~$\bar{x}$}.
  When $X$ is a model over~$\Z_{(\ell)}$ of a projective variety over~$\Q$,
  then $X(\Z_{(\ell)}) = X(\Q)$, so that we can think of the residue class
  of~$\bar{x}$ as the set of rational points on~$X$ reducing mod~$\ell$
  to~$\bar{x}$.
\end{definition}

Recall that $X_1(p)$ has a smooth model over~$\Z[\frac{1}{p}]$; this implies
the corresponding statement for the $d$th symmetric power~$X_1(p)^{(d)}$.

\begin{lemma} \label{lem:crit_residue_classes}
  Let $\ell \neq p$ be a prime. Assume that
  \begin{enumerate}[\upshape(a)]
    \item \label{rc_a}
          The residue class of each point $\bar{x} \in X_1(p)^{(d)}(\F_\ell)$
          that is a sum of images under~$\red_\ell$ of rational cusps
          contains at most one rational point.
    \item \label{rc_b}
          The residue class of each point $\bar{x} \in X_1(p)^{(d)}(\F_\ell)$
          that is not a sum of images under~$\red_\ell$ of rational cusps
          contains no rational point.
  \end{enumerate}
  Then $p \notin S(d)$.
\end{lemma}

\begin{proof}
  Let $x \in X_1(p)^{(d)}(\Q)$ and write
  $\bar{x} = \red_\ell(x) \in X_1(p)^{(d)}(\F_\ell)$.
  By assumption~\eqref{rc_b}, $\bar{x} = \bar{x}_1 + \dots + \bar{x}_d$
  is a sum of images of rational cusps. Let $x_1, \ldots, x_d \in X_1(p)(\Q)$
  be rational cusps such that $\red_\ell(x_j) = \bar{x}_j$ for
  $1 \le j \le d$. Then $x' = x_1 + \dots + x_d \in X_1(p)^{(d)}(\Q)$
  is such that $\red_\ell(x') = \bar{x}$. By assumption~\eqref{rc_a},
  $x$ is the only rational point in the residue class of~$\bar{x}$,
  so it follows that $x' = x \in \im(\alpha)$ with $\alpha$ as
  in~\cref{lem:criterion_alpha}. So $\alpha$ is surjective, and
  \cref{lem:criterion_alpha} shows that $p \notin S(d)$.
\end{proof}

Fix a rational cusp $c \in X_1(p)(\Q)$. We can then define a morphism
$\iota \colon X_1(p)^{(d)} \to J_1(p)$ by sending $x_1 + \dots + x_d$
to the class of the divisor $x_1 + \dots + x_d - d \cdot c$; here $J_1(p)$
denotes the Jacobian variety of~$X_1(p)$; see~\cref{sec:mod_curves} below.
This map is actually defined over~$\Z[\frac{1}{p}]$.

The standard way of verifying assumption~\eqref{rc_a} is to
show that there is a morphism of abelian varieties $t \colon J_1(p) \to A$
such that
\begin{enumerate}[(i)]
  \item $t \circ \iota$ is injective on each residue class of
        a point~$\bar{x}$ as in assumption~\eqref{rc_a}, and
  \item $\red_\ell \colon t(J_1(p)(\Q)) \to A(\F_\ell)$ is injective.
\end{enumerate}
By standard properties of $\red_\ell$ on the rational torsion subgroup,
the second condition is satisfied when $t(J_1(p)(\Q))$ is finite
and either $\ell$ is odd or $\ell = 2$ and $t(J_1(p)(\Q))$ has odd order.
We can achieve this by choosing $A$ as a factor of~$J_1(p)$ that has
Mordell-Weil rank zero and $t$ to be the projection to~$A$ (plus some
technicalities when $\ell = 2$). By work of
Kolyvagin-Logach\"ev~\cite{kolyvagin-logachev} and Kato~\cite{kato},
it is known that the ``winding quotient''~$J^{\bfe}_1(p)$ of~$J_1(p)$
has Mordell-Weil rank zero. Assuming the Birch and Swinnerton-Dyer Conjecture
for abelian varieties, $J^{\bfe}_1(p)$ is in fact the largest such quotient.
See~\cref{sec:mod_curves} for the definition of the winding quotient.

The first condition follows if it can be shown that $t \circ \iota$
is a ``formal immersion'' at the relevant points~$\bar{x}$;
see~\cref{sec:formal_imm}.

We can work with~$J_0(p)$ in place of~$J_1(p)$. Then there is
only one point~$\bar{x}$ to consider, which is $d$~times the image
of the rational cusp~$\infty$ on~$X_0(p)$. This is what Mazur and Kamienny
used to determine $S(1)$ and~$S(2)$ and is also used in Merel's proof
of an explicit bound on~$S(d)$ for all~$d$ and Oesterl\'e's improvement
of the bound. In all this work, odd primes~$\ell$ are used.
To deal with~$S(3)$, Parent had to work with~$J_1(p)$ (which was made
possible by Kato's work showing that the winding quotient has rank zero)
and also had to use $\ell = 2$ to exclude some of the primes.

One minor innovation we introduce here is that we work with some
intermediate curve~$X_H$ between $X_1(p)$ and~$X_0(p)$; see
again~\cref{sec:mod_curves}. This can reduce the necessary work in
cases when using~$J_0(p)$ is not successful, but the dimension of~$J_1(p)$
is too large to make computations feasible.

Assuming Oesterl\'e's bound~\eqref{eq:oesterle}, verification of
assumption~\eqref{rc_a} amounts to exhibiting a suitable~$t$ for
each prime $p \le (3^{d/2} + 1)^2$ such that $p \notin S(d)$ and
checking that it satisfies the conditions. This can be done by
an explicit computation using modular symbols, which is based on
a criterion established by Kamienny for~$J_0(p)$ and extended to~$J_1(p)$
by Parent. In view of assumption~\eqref{rc_b} (see below), we
work with $\ell = 2$, which necessitates using ``Parent's trick''
to deal with the technicalities that arise when $\ell$ is not odd.

For certain small primes, this is not sufficient. For $d \le 7$,
these primes~$p$ have the property that $J_1(p)(\Q)$ is finite,
which allows us to work with the full Jacobian and perform some
more direct computations. This is another new ingredient compared to
earlier work. In the course of our work, we establish an open
case of a conjecture of Conrad, Edixhoven and Stein: we show that
the group~$J_1(29)(\Q)$ (which is finite) is generated by differences
of rational cusps on~$X_1(29)$; see~\cref{thm:29cusp}.

Combining both approaches, we obtain the following result.

\begin{proposition} \label{prop:ass_a}
  Let $p \le 2281 = \floor{(3^{7/2} + 1)^2}$ be a prime. If
  \[ \begin{array}{*{2}{r@{\quad\text{and}\quad}l@{\qquad\text{or}\qquad}}r@{\quad\text{and}\quad}l}
        d = 3 & p \ge 17 & d = 4 & p \ge 19 & d = 5 & p \ge 23 \\[6pt]
        \multicolumn{2}{r@{\qquad}}{\text{or}} & d = 6 & p \ge 23 & d = 7 & p \ge 29,
     \end{array}
  \]
  then assumption~\eqref{rc_a} of~\cref{lem:crit_residue_classes}
  is satisfied.
\end{proposition}

\cref{prop:ass_a} is proved in~\cref{sec:ass_a}.

We now consider assumption~\eqref{rc_b} of~\cref{lem:crit_residue_classes}.
The simplest way for the assumption to be satisfied is when there are no
points~$\bar{x}$ that are not sums of images of rational cusps. Equivalently,
\begin{enumerate}[(i)]
  \item \label{alpha_cond1}
        there is no elliptic curve~$E$ over $\F_{\ell^{d'}}$ with $d' \le d$
        such that $p \mid \#E(\F_{\ell^{d'}})$, and
  \item \label{alpha_cond2}
        $p \nmid \ell^{d'} \pm 1$ for all $d' \le d$.
\end{enumerate}
The first condition excludes the existence of non-cuspidal points,
whereas the second excludes the possibility that $X_1(p)(\F_{\ell^{d'}})$
contains cusps that are not images of rational cusps. Recall that the
irrational cusps are defined over the maximal real subfield of~$\Q(\mu_p)$,
which has a place of degree dividing~$d$ above~$\ell$ if and only if
$\ell^d \equiv \pm 1 \bmod p$.

We note the following simple consequence.

\begin{lemma} \label{lem:alpha_surj_simple}
  If $p > (\ell^{d/2} + 1)^2$, then assumption~\eqref{rc_b}
  of~\cref{lem:crit_residue_classes} is satisfied.
\end{lemma}

\begin{proof}
  If there is an elliptic curve~$E$ over $\F_{\ell^{d'}}$ with $d' \le d$
  such that $p \mid \#E(\F_{\ell^{d'}})$, then by the Hasse bound,
  \[ p \le \#E(\F_{\ell^{d'}}) \le (\ell^{d'/2} + 1)^2 \le (\ell^{d/2} + 1)^2 , \]
  which is not the case, so condition~\eqref{alpha_cond1} above is satisfied.
  Since $p > (\ell^{d/2} + 1)^2 > \ell^d + 1$, condition~\eqref{alpha_cond2}
  is also satisfied.
\end{proof}

This explains the form of Oesterl\'e's bound~\eqref{eq:oesterle}, which
is related to the fact that he is working with $\ell = 3$.

We also see that it is advantageous to use the smallest possible~$\ell$,
because then the condition of~\cref{lem:alpha_surj_simple} covers more
primes~$p$. But even using $\ell = 2$, we need to verify
assumption~\eqref{rc_b} for some primes $p < (2^{d/2} + 1)^2$.
In some cases, we can still show for such primes that there are no
points~$\bar{x}$ that are not sums of images of rational cusps,
but this is not enough: when
\begin{align*}
  (d, p) \in \{&(5, 31), (5, 41), (6, 29), (6, 31), (6, 41), (6, 73), \\
               &(7, 29), (7, 31), (7, 37), (7, 41), (7, 43), (7, 59), \\
               &(7, 61), (7, 67), (7, 71), (7, 73), (7, 113), (7, 127)\} ,
\end{align*}
there actually are such points, and we have to work quite a bit harder
to show that they are not images of rational points on~$X_1(p)^{(d)}$.
This is another novel aspect of our work. We use a number of different
approaches (for $p = 37$, see further below).
\begin{enumerate}[(1)]\addtolength{\itemsep}{3pt}
  \item \label{b_item_1}
        For $p \in \set{29, 31, 41}$, we can again use direct computations
        based on the fact that $J_1(p)(\Q)$ is finite and known;
        see~\cref{lem:29_31_41_alpha}
  \item \label{b_item_2}
        For $p \in \set{71, 113, 127}$ and $d = 7$, we use a new criterion
        based on gonality estimates and working with Hecke operators
        as correspondences, which shows directly that $p \notin S(d)$;
        see~\cref{cor:7_71_113_127}.
  \item \label{b_item_3}
        For $(d, p) \in \set{(6, 73), (7, 43)}$, we use an intermediate
        curve~$X_H$ such that $X_H^{(d)}$ possesses a rational point~$x_H$
        in the image of the relevant residue class and use a formal immersion
        argument to show that it is the only rational point in this residue
        class. This implies that every rational point on~$X_1(p)^{(d)}$ in
        the residue class of~$\bar{x}$ must map to~$x_H$, but $x_H$ does not
        lift to a rational point on~$X_1(p)^{(d)}$;
        see~\cref{lem:6_73,lem:7_43}.
  \item \label{b_item_4}
        For $p \in \set{59, 61, 67, 73}$ and $d = 7$, we use another new
        criterion that shows that a non-cuspidal point
        $\bar{x} \in X_1(p)^{(d)}(\F_2)$ is not the reduction mod~$2$
        of a rational point by showing that its image in~$J_1(p)(\F_2)$
        is not in the reduction of the Mordell-Weil group; see~\cref{lem:7_p}.
\end{enumerate}

We then obtain the following result.

\begin{proposition} \label{prop:ass_b}
  For the following pairs of an integer $3 \le d \le 7$ and a prime~$p$,
  assumption~\eqref{rc_b} of~\cref{lem:crit_residue_classes} is satisfied.
  \begin{align*}
    d = 3 \colon& \qquad p = 11 \quad\text{or}\quad p \ge 17 \\
    d = 4 \colon& \qquad p \ge 19 \\
    d = 5 \colon& \qquad p \ge 23 \\
    d = 6 \colon& \qquad p \ge 23 \quad\text{and}\quad p \neq 37 \\
    d = 7 \colon& \qquad p \ge 29 \quad\text{and}\quad p \neq 37
  \end{align*}
\end{proposition}

\cref{prop:ass_b} is proved in~\cref{sec:ass_b}.

We still have to show~\cref{prop:6_37} and that $37 \notin S(7)$.
We combine the approaches in \eqref{b_item_3} and~\eqref{b_item_4}
to do this.
We first show using~\eqref{b_item_4} that no non-cuspidal point
in~$X_1(37)^{(7)}(\F_2)$ is the reduction mod~$2$ of a rational
point and there is essentially only one such point
in~$X_1(37)^{(6)}(\F_2)$. We then use the formal immersion argument
as in~\eqref{b_item_3} to show that the remaining points
in~$X_1(37)^{(d)}(\F_2)$ for $d = 6, 7$ lift uniquely to rational
points; see~\cref{lem:37_6,lem:37_7}.

\cref{thm:main} then follows from this and~\cref{prop:ass_a,prop:ass_b}
using~\cref{lem:crit_residue_classes} and Oesterl\'e's
bound~\eqref{eq:oesterle}.

A large part of the work done in this paper relies heavily on computations.
We provide Magma~\cite{Magma} code (with explanatory comments)
for all these computations at~\cite{code}.
The timings we give in some places in this paper were obtained on the
last authors's current laptop (as of~2020). All computations together
took about one day on this machine. We also provide SageMath~\cite{sage} code
at the first author's Github site~\cite{sage-code} that independently
verifies the claims made in~\cref{sec:ass_a}.
Some of these computations rely on modular symbols.
See for example~\cite{stein} for the necessary background.


\subsection*{The structure of the paper}

We begin by recalling some background on modular curves
in~\cref{sec:mod_curves}. In~\cref{sec:small}, we quote the list
of primes~$p$ such that $J_1(p)(\Q)$ is finite from~\cite{CES}
and prove that for such primes, $J_1(p)(\Q)$ is generated by
differences of rational cusps (the new case being $p = 29$)
and that the reduction map $J_1(p)(\Q) \to J_1(p)(\F_2)$ is injective.
We use this to prove assumption~\eqref{rc_b} for $p = 29$, $31$, and~$41$.
In~\cref{sec:formal_imm}, we introduce formal immersions
and state the computational criterion we use to verify
assumption~\eqref{rc_a}. \cref{sec:ass_a} reports on these computations,
and~\cref{sec:oesterle} contains the proof of Oesterl\'e's
bound~\eqref{eq:oesterle}. In~\cref{sec:modlarge} we state and prove
the criterion used to show that $71, 113, 127 \notin S(7)$.
Finally, we complete the verification
of assumption~\eqref{rc_b} in~\cref{sec:ass_b}, which also
contains the proof of~\cref{prop:6_37}.


\subsection*{What is new in this paper?}

The main new \emph{result} is~\cref{thm:main}, which extends the list
of known sets~$S(d)$ from $d \le 3$ to $d \le 7$.
Completing the determination of~$S(6)$, \cref{prop:6_37} gives a
classification of the sporadic points in~$X_1(37)^{(6)}(\Q)$.
Another new result
is~\cref{thm:29cusp}, which confirms a conjecture made in~\cite{CES}
in a case that was left open in that paper.

We also develop some new \emph{techniques} for proving that $p \notin S(d)$
for suitable $d \ge 1$ and primes~$p$. One point is the use of intermediate
curves in various computations instead of just either~$X_0(p)$ or~$X_1(p)$.
Another is the use of explicit computations in the Picard group
of~$X_1(p)_{\F_2}$ when $p$ is such that $J_1(p)(\Q)$ is finite.
In addition, we derive two new criteria, one that uses the gonality of~$X_1(p)$
and can show directly that $p \notin S(d)$ using global arguments
(\cref{prop:hecke_gon}),
and a related one that works over~$\F_2$ using Hecke correspondences
(\cref{lem:filter}).
Finally, we extend the formal immersion approach that is traditionally
used to show what we call assumption~\eqref{rc_a} to also apply to
assumption~\eqref{rc_b}. All this is necessary to be able to
determine~$S(7)$.


\subsection*{Why stop at $d = 7$?}

Obviously, determining~$S(d)$ gets harder and harder as $d$ grows.
When $d = 1$, the formal immersion condition for $X_0(p)$ is essentially
trivially satisfied, and assumption~\eqref{rc_b} for $\ell = 3$ is
automatically satisfied for
$p > \bigl\lfloor{(\sqrt{3} + 1)^2}\bigr\rfloor = 7$.
Once the theoretical framework is in place (which, of course, was Mazur's
key contribution in~\cites{mazur1,mazur2}), no computation is necessary
to obtain the desired result.

For $d = 2$, Kamienny had to come up with a criterion that allows to
verify the formal immersion condition (still for~$X_0(p)$).
In this case, it can still be
shown to hold by a theoretical argument. The trivial bound for
assumption~\eqref{rc_b} when $\ell = 3$ is $p > 16$, which is again
sufficient.

For $d = 3$ and larger, one needs to work to verify the formal immersion
condition. Merel and Oesterl\'e managed to find a theoretical argument that
does this (for $X_0(p)$ and $\ell = 3$) for $p$ larger than some explicit
polynomial in~$d$. Oesterl\'e then came up with another ingenious way to
reduce the remaining cases to a finite and manageable amount of computation,
thus proving the bound~\eqref{eq:oesterle}. To determine~$S(3)$, Parent
had to rely on this and to come up with a way of using~$X_1(p)$ and
$\ell = 2$ to cover the primes between $\ceil{(2^{3/2} + 1)^2} = 15$
and $\floor{(3^{3/2} + 1)^2} = 38$ (and $p = 43$, which had escaped
Oesterl\'e's approach).

For $d \ge 4$, there are two main difficulties that each get worse
as $d$ increases.
\begin{enumerate}[(1)]\addtolength{\itemsep}{3pt}
  \item The gap between the best general bound~\eqref{eq:oesterle}
        and the smallest prime not in~$S(d)$ increases
        exponentially with~$d$. While we can, for each~$d$ and
        each~$p$ in this range, verify the formal immersion condition
        for~$X_0(p)$ or some intermediate curve~$X_H$ computationally,
        the computational effort increases considerably with~$p$.
        For $d = 7$, this part of the computation took about two hours.
        For $d = 8$, the upper end of this range is about three times
        as large as for $d = 7$, which lets us expect that doing this in
        reasonable time would require a massively parallel computation.
        For $d \ge 9$, this appears to be infeasible in the absence
        of a major theoretical advance that leads to a significantly
        reduced general bound.
  \item The gap between the primes in~$S(d)$ and the ``easy'' range
        for assumption~\eqref{rc_b} also increases. Most likely, this
        increase is also exponential, since we expect that $\max S(d)$
        should grow only polynomially (possibly even linearly).
        This means that there will be more and more primes~$p$ for
        which we have to show assumption~\eqref{rc_b} when there
        are indeed points in~$X_1(p)^{(d)}(\F_2)$ that are not sums
        of images of rational cusps. While we could deal with the
        ``rank zero primes'' $p = 29, 31, 41$ by explicit computations
        and with the one further such prime $p = 73$ for $d = 6$
        by a variant of the formal immersion criterion,
        this is the point where it gets hard when $d = 7$. To rule out
        the primes $37$, $43$, $59$, $61$, $67$, $71$, $73$, $113$,
        and~$127$, we needed to develop some new criteria, and some
        of the computations that are then still necessary run for
        several hours.
\end{enumerate}

However, it appears that our new criteria can be used to go a bit
further. This will be explored in a follow-up paper.


\subsection*{Acknowledgments}

We would like to thank Bas Edixhoven, Barry Mazur, and Lo\"ic Merel for their
many valuable comments and suggestions, Pierre Parent for his idea to look at
CM elliptic curves for the proof of $73 \notin S(6)$ (\cref{lem:6_73}),
Filip Najman for some helpful
information on sporadic torsion points, and Tessa Schild for her proofreading
of an earlier version of this paper. We thank Joseph Oesterl\'e for kindly
allowing us to use his notes~\cite{oesterle} and helping the first author
understand the proof of the bound~\eqref{eq:oesterle}.
We also thank the anonymous
referee of an earlier version of this paper for some valuable feedback.


\section{Preliminaries on modular curves} \label{sec:mod_curves}

A good reference for most of the following is~\cite{diamondim}.

As usual, we define, for $N \in \Z_{\ge 1}$,
\begin{align*}
  \Gamma_0(N) &= \left\{\smm{a}{b}{c}{d} \in \SL_2(\Z) : N \mid c\right\}
   \quad\text{and} \\
  \Gamma_1(N) &= \left\{\smm{a}{b}{c}{d} \in \SL_2(\Z) : (c, d) \equiv (0, 1) \bmod N\right\} .
\end{align*}
$\SL_2(\Z)$ and therefore also $\Gamma_0$ and~$\Gamma_1$ act on the
complex upper half plane~$\HH$ and on $\HH^* = \HH \cup \P^1(\Q)$ by
M\"obius transformations.
Then the quotient $Y_j(N)(\C) = \Gamma_j(N) \backslash \HH$ (for $j = 0, 1$)
is a Riemann surface that can be compactified to
$X_j(N)(\C) = \Gamma_j(N) \backslash \HH^*$ by adding the finitely many
cusps $C_j(N) = \Gamma_j(N) \backslash \P^1(\Q)$. The points in~$Y_1(N)(\C)$
classify pairs~$(E, P)$ consisting of an elliptic curve~$E$ over~$\C$
and a point $P \in E(\C)$ of exact order~$N$; in terms of a representative
point $\tau \in \HH$, this is given by $E = \C/\Lambda_\tau$ with
$\Lambda_\tau = \Z + \Z \tau$ and $P = \tfrac{1}{N} + \Lambda$.
Similarly, $Y_0(N)(\C)$ classifies pairs~$(E, C)$ where again $E$ is an
elliptic curve over~$\C$ and $C \subseteq E(\C)$ is a cyclic subgroup
of order~$N$. The compact Riemann surfaces~$X_j(N)(\C)$ can be identified
with the set of complex points on projective algebraic curves~$X_j(N)$
defined over~$\Q$ (or even over~$\Z[\tfrac{1}{N}]$). The rational structure
is defined in terms of the $q$-expansions of functions on~$X_j(N)(\C)$:
such a function~$f$ lifts to a modular function with respect to~$\Gamma_j(N)$
on~$\HH$ and therefore has a Laurent series expansion in terms of
$q = e^{2\pi i \tau}$. The function field $\Q(X_j(N))$ is then defined to
consist of those~$f$ whose $q$-expansion has coefficients in~$\Q$.
Since the rational structure is defined in terms of~$q$, the natural
moduli interpretation of a point on~$X_1(N)$ over~$\Q$ (or any field
extension~$K$) is as representing a pair~$(E, \varphi)$, where $E$
is an elliptic curve over~$\Q$ (or~$K$) and $\varphi \colon \mu_N \to E$
is an embedding of the group~$\mu_N$ of $N$th roots of unity into~$E$
as group schemes. This is because the image of~$\tfrac{1}{N}$ under
$\tau \mapsto e^{2 \pi i \tau}$ is not rational, but a generator of~$\mu_N$.
Since (over any field~$K$ of characteristic not dividing~$N$) there is a
natural bijection between pairs $(E, P)$ and pairs $(E', \varphi)$ as above,
the points on~$X_1(N)$ can still be understood as classifying elliptic
curves over~$K$ together with a point of order~$N$, but we have to keep in
mind that this is not the same as the moduli interpretation over~$\C$
given above. (The bijection is obtained as follows. Given a pair~$(E, P)$
and $\zeta \in \mu_N$, the set of $Q \in E[N]$ such that $e_N(Q, P) = \zeta$
forms a coset $C_\zeta$ of the subgroup~$\Z P$ generated by~$P$. We then set
$E' = E/\Z P$ and $\varphi \colon \zeta \mapsto C_\zeta/\Z P$.)

The space of cusp forms for~$\Gamma_j(N)$ is canonically isomorphic
to the space of regular differentials on~$X_j(N)(\C)$. Under this
isomorphism, regular differentials on~$X_j(N)$ over~$\Q$ correspond
to cusp forms whose $q$-expansion has rational coefficients.

There is a natural map $X_1(N) \to X_0(N)$ (induced over~$\C$ by the
identity on~$\HH^*$). This makes $X_1(N)$ into a (possibly ramified) Galois
covering of~$X_0(N)$, whose Galois group consists of the diamond operators
$\diamondop{a}$ for $a \in (\Z/N\Z)^\times$, where $\diamondop{-1}$
is the identity, so the Galois group is naturally isomorphic
to $(\Z/N\Z)^\times/\set{\pm 1}$. In terms of the interpretation
of points on~$Y_1(N)$ as pairs $(E, \varphi \colon \mu_N \to E)$,
the action of~$\diamondop{a}$ corresponds to pre-composing $\varphi$
with the $a$th-power map.
If $H \subseteq (\Z/N\Z)^\times/\set{\pm 1}$
is a subgroup, then we have an intermediate curve $X_H = H \backslash X_1(N)$
between $X_1(N)$ and~$X_0(N)$.

We write $\infty \in X_j(N)$ for the cusp that over~$\C$ is the image
of $\infty \in \P^1(\Q)$. Note that $\infty \in X_j(N)(\Q)$, since
it corresponds to $q = 0$. When $N = p$ is prime, $X_0(p)$ has the two
cusps~$\infty$ and the cusp represented by $0 \in \P^1(\Q)$, which
are both rational, whereas $X_1(p)$ has $p-1$ cusps, which split into
two orbits under the diamond operators, each consisting of $(p-1)/2$
cusps. One of the orbits contains~$\infty$ and consists of rational
cusps, the other orbit consists of cusps defined over the maximal
totally real subfield of the cyclotomic field~$\Q(\mu_p)$; these
cusps are all conjugate under the Galois action, and the Galois action
is given by diamond operators (since it commutes with them).
An analogous statement is true for the cusps of~$X_H$.
See~\cite{stevens}*{Thm.~1.3.1} for a description of the Galois
action on the cusps.

We denote the Jacobian varieties of $X_0(N)$, $X_1(N)$ and~$X_H$
by $J_0(N)$, $J_1(N)$ and~$J_H$, respectively. They are defined
over~$\Q$ and extend to abelian schemes over~$\Z[\tfrac{1}{N}]$.

We denote the Hecke algebra, in its various incarnations, by~$\T$.
It is generated by the Hecke operators~$T_n$ for all $n \ge 1$,
or alternatively, by all $T_p$ for $p$ prime together with
diamond operators~$\diamondop{a}$ for $a$ generating
$(\Z/N\Z)^\times/\{\pm 1\}$. The Hecke algebra acts on the
integral homology $H_1(X_H(\C), \Z)$, the relative homology
$H_1(X_H(\C), \cusps, \Z)$, the associated spaces of modular forms
or cusp forms, and as endomorphisms of~$J_H$. The Hecke operators~$T_n$
and the diamond operators~$\diamondop{a}$ can also be viewed as
correspondences on~$X_H$. It will
always be clear from the context or explicitly stated which
interpretation is considered.

The integral relative homology with respect to the cusps is
generated as a $\Z$-module by modular symbols $\set{\gamma_1, \gamma_2}$
with $\gamma_1, \gamma_2 \in \P^1(\Q)$. There is an integration
pairing
\[ H_1(X_H(\C), \cusps, \Z) \times H^0(X_{H,\C}, \Omega^1) \to \C, \qquad
   (\xi, \omega) \mapsto \int_\xi \omega
\]
(if $\xi = \set{\gamma_1, \gamma_2}$, then the integral is along any path
in~$\HH^*$ joining $\gamma_1$ to~$\gamma_2$);
it induces a perfect pairing of real vector spaces between
$H_1(X_H(\C), \R)$ and $H^0(X_{H,\C}, \Omega^1)$, and the composition
\[ \pi \colon H_1(X_H(\C), \cusps, \Z) \to H^0(X_{H,\C}, \Omega^1)^* \to H_1(X_H(\C), \R) \]
has image in the rational homology~$H_1(X_H(\C), \Q)$ by the
Manin-Drinfeld Theorem \cites{manin,drinfeld}.

\begin{definition} \label{def:winding}
  We set
  \[ \bfe = \pi\bigl(-\set{0, \infty}\bigr) \in H_1(X_H(\C), \Q) ; \]
  this is called the \emph{winding element}.
  Its annihilator~$\Ann(\bfe)$ in~$\T$ is the \emph{winding ideal}.
  It acts via endomorphisms on~$J_H$; the quotient
  $J^{\bfe}_H := J_H/\Ann(\bfe) J_H$ is the \emph{winding quotient}.
\end{definition}

The defnition of the winding element goes back to
Mazur~\cite{mazur1}*{Lemma~II.18.6 and the definition preceding it}
in the case of~$J_0(N)$.
We note that there is some ambiguity regarding the sign of the winding
element in the literature. We follow~\cite{merel}*{Section~1} here
(but, for example, \cite{parent1} uses the opposite sign.)
The winding quotient has the following essential property.

\begin{theorem} \label{thm:winding}
  For each subgroup $H \subseteq (\Z/N\Z)^\times/\set{\pm 1}$, the
  Mordell-Weil group $J^{\bfe}_H(\Q)$ is finite.
\end{theorem}

Merel \cite{merel}*{\S1} was the first one to introduce the winding
quotient for~$J_0(p)$ with $p$~prime, where he also proves that its
Mordell-Weil group is finite using a result from~\cite{kolyvagin-logachev},
which states that an abelian variety~$A$ over~$\Q$ that is a quotient
of~$J_0(N)$ has Mordell-Weil rank~$0$ when $L(A, 1) \neq 0$.
Parent in~\cite{parent1}*{\S3.8} generalized Merel's statement to composite
numbers~$N$. The result of Kolyvagin and Logach\"ev was generalized by
Kato~\cite{kato}*{Cor.~14.3} to quotients
of~$J_1(N)$. In both \cite{parent2} and~\cite{parent3}, it is mentioned
that the theorem follows from Kato's generalization. This can be seen by
adapting the arguments of~\cite{parent1}*{\S3.8} accordingly.
The key point in the proof is that $J^{\bfe}_H$ is isogenous to a product
of simple abelian varieties~$A$ over~$\Q$ such that $L(A, 1) \neq 0$.
Kato's result then shows that $A(\Q)$ is finite.

The following is a variant of~\cite{parent2}*{Prop.~1.8}. We remark
that, according to~\cite{diamondim}*{p.~87}, the Eichler-Shimura
relation on~$X_1(N)$ with the modular interpretation used here
(and in~\cite{parent2}) is different from that valid with the more
usual interpretation as parameterizing pairs~$(E, P)$ of elliptic
curves with a point of order~$N$. We therefore believe that our
version is correct, and that (the first part of) Parent's statement
needs to be modified accordingly.

\begin{proposition} \label{prop:ann_rat_tors}
  Let $q \nmid N$ be a prime and $P \in J_H(\Q)_{\tors}$ such that
  $q$ is odd or $P$ has odd order. Then $(T_q - \diamondop{q} - q)(P) = 0$.
\end{proposition}

\begin{proof}
  Let $n$ be the order of~$P$.
  Then $(T_q - \diamondop{q} - q)(P) \in J_H(\Q)$ is a point of order dividing~$n$.
  We write $\bar{P}$ for the reduction mod~$q$ of~$P$,
  $\Frob_q$ for the Frobenius on~$J_{H,\F_q}$ and $\Ver_q$ for its dual
  (Verschiebung). Then we have the Eichler-Shimura relation
  \[ T_{q,\F_q} = \diamondop{q} \Frob_q + \Ver_q , \qquad\text{and}\qquad
     \Ver_q \circ \Frob_q = q
  \]
  in $\End_{\F_q}(J_{H,\F_q})$; see~\cite{diamondim}*{p.~87}. So,
  using that $\Frob_q(\bar{P}) = \bar{P}$,
  \[ T_{q,\F_q}(\bar{P}) = \diamondop{q} \Frob_q(\bar{P}) + \Ver_q(\bar{P})
                         = \diamondop{q} \bar{P} + q \bar{P} ,
  \]
  which implies that $(T_{q,\F_q} - \diamondop{q} - q)(\bar{P}) = 0$.
  Since the reduction map is injective on~$J_H(\Q)_{\tors}$ when
  $q$ is odd, and it is injective on odd order torsion when $q = 2$,
  the claim follows.
\end{proof}

\begin{remark}
  We note that for the proof of~\cref{thm:29cusp}, it actually does
  not matter whether one uses $T_q - \diamondop{q} - q$ as
  in~\cref{prop:ann_rat_tors} or $T_q - \diamondop{q} q - 1$
  as in~\cite{parent2}. Up to composition with~$\diamondop{q}$
  or its inverse, the two operators are conjugate to each other
  under the Atkin-Lehner involution~$w_p$;
  see~\cite{diamondim}*{p.~56 and Rmk.~10.2.2}.
  If we use the ``wrong'' operators in the proof of~\cref{thm:29cusp},
  then instead of $J_1(29)(\Q)_{\tors} \subseteq C$,
  we find that $w_{29}(J_1(29)(\Q)_{\tors}) \subseteq C$, which
  also implies that
  $J_1(29)(\Q)_{\tors} \subseteq w_{29}(C) = C$ (the cusps are permuted
  by~$w_p$).

  Our second application in~\cref{cor:trick}
  uses the operators with $q$ odd to kill torsion. By the remark
  after~\cref{cor:beta_inj}, it is enough to kill $2$-torsion.
  For $q$ an odd prime, the two operators $T_q - \diamondop{q} - q$
  and $T_q - \diamondop{q} q - 1$ differ by a multiple of~$2$ in the
  Hecke algebra, so they have the same effect on $2$-torsion points.
  This implies that the conclusion of~\cref{cor:beta_inj} also holds
  if we use the ``wrong'' operator and show that $t \circ \iota$
  is a formal immersion. In particular, the conclusions of~\cite{parent2}
  are valid.
\end{remark}

We will also need the following statement.

\begin{proposition}[Derickx] \label{prop:hecke_op_kernel}
  Let $q \nmid N$ be a prime. We consider $t = T_q - \diamondop{q} - q \in \T$
  as a correspondence on~$X_1(N)$, inducing an endomorphism of the
  divisor group of~$X_1(N)$ over~$\C$. Then the kernel of~$t$ is
  contained in the subgroup of divisors supported in cusps.
\end{proposition}

\begin{proof}
  Let $D$ be a divisor in the kernel of~$t$, so that
  \begin{equation} \label{eq:hecke_rel}
    T_q(D) = \diamondop{q}(D) + q D \,.
  \end{equation}
  A non-cuspidal point $x \in X_1(N)(\C)$ corresponds to an elliptic
  curve~$E$ over~$\C$ with additional structure.
  The point $\diamondop{q}(x)$ corresponds to
  the same curve~$E$ (with modified extra structure), and $T_q(x)$ is a
  sum of points corresponding to all the elliptic curves that are
  $q$-isogenous to~$E$. We define the \emph{$q$-isogeny graph}~$G$
  to have as vertices the isomorphism classes of all elliptic curves
  over~$\C$; two vertices are connected by an edge when there is a
  $q$-isogeny between the corresponding curves.
  There is a natural map $\gamma$ from $Y_1(N)(\C)$
  to the vertex set of~$G$. Let $x$ be a non-cuspidal point in the support
  of~$D$ and let $G_x$ be the connected component of~$G$
  containing~$\gamma(x)$. Let $E$ be the elliptic curve given by~$x$.
  We distinguish two cases.

  First, assume that $E$ does not have~CM. Then $G_x$ is an infinite
  $(q+1)$-regular tree. Consider a vertex~$v$ of~$G_x$ that has maximal
  possible distance from~$\gamma(x)$ among all vertices of the
  form~$\gamma(y)$ for a point~$y$ in the support of~$D$.
  Let $y_1, \ldots, y_n$ be the points in the support of~$D$ such that
  $\gamma(y_j) = v$, and let $w$ be a neighbor of~$v$ whose distance
  from~$\gamma(x)$ is larger than that of~$v$. Each $T_q(y_j)$ contains
  precisely one point~$y'_j$ such that $\gamma(y'_j) = w$, and these
  points are distinct for distinct points~$y_j$. Since $w$ is not of
  the form~$\gamma(z)$ for a point~$z$ in the support of~$D$, this shows
  that $T_q(D)$ has points in its support that do not occur in the
  support of $\diamondop{q}(D) + q D$ (recall that
  \hbox{$\gamma(\diamondop{q}(y)) = \gamma(y)$}). This contradicts the
  relation~\eqref{eq:hecke_rel}, and we conclude that there can be no
  non-CM point~$x$ in the support of~$D$.

  Now consider the case that $E$ has~CM. Then $G_x$ is no longer a
  tree in general, but has the structure of a ``volcano'';
  see~\cite{sutherland-ants}. For a CM elliptic curve over~$\C$,
  this volcano has infinite depth. Concretely, this means that it
  consists of a number of rooted $(q+1)$-regular trees whose roots
  form a cycle (which may have length $1$ or~$2$). We can now argue
  as in the first case by choosing~$v$ to be a vertex of maximal
  level (i.e., distance from the root cycle) and~$w$ to be a neighbor
  of~$v$ whose level is larger by one. This shows that there can be
  no CM points in the support of~$D$ as well.

  The only points that we have not excluded from the support of~$D$
  are the cusps; this proves the claim.
\end{proof}

\begin{remark}
  In the case that $N = p$ is a prime, we can describe the kernel exactly.
  The rational cusps are killed by~$t$, whereas the irrational cusps
  are killed by $t^* = T_q - q \diamondop{q} - 1$;
  compare~\cite{parent2}*{Section~2.4} (the rational cusps are
  those mapping to the cusp~$\infty$ on~$X_0(p)$ under the modular
  interpretation we use). Since $t - t^* = (q-1)(\diamondop{q} - 1)$
  and the divisor group is torsion-free, $t$ kills a divisor supported
  on irrational cusps if and only if it is invariant under~$\diamondop{q}$.
\end{remark}


\section{Rank zero primes} \label{sec:small}

We say that a prime~$p$ is a \emph{rank zero prime}
when $J_1(p)(\Q)$ is finite.

The following result gives us the list of rank zero primes.
This is~\cite{CES}*{Prop.~6.2.1}; we include some more information from
Section~6.2 of~loc.~cit.

\begin{proposition} \label{prop:rank0}
  The rank zero primes~$p$ are the primes
  $p \le 31$ and $41$, $47$, $59$, and~$71$.

  For all of these, except possibly $p = 29$, the group~$J_1(p)(\Q)$
  is generated by differences of rational cusps, and for all except
  $p = 17$, $29$, $31$ and~$41$, the order of~$J_1(p)(\Q)$ is odd.
\end{proposition}

We can add to this the following new result, which confirms
Conjecture~6.2.2 in~\cite{CES} for the smallest open case $p = 29$.

\begin{theorem} \label{thm:29cusp}
  The group $J_1(29)(\Q)$ is generated by differences of rational cusps.
\end{theorem}

\begin{proof}
  We prove this by a computation using modular symbols, as follows.
  The group $J_1(29)(\C)_{\tors}$ is canonically isomorphic to
  $M := H_1(X_1(29)(\C), \Z) \otimes_{\Z} \Q/\Z$. By~\cref{prop:ann_rat_tors},
  the image of the rational torsion subgroup is annihilated
  by $T_q - \diamondop{q} - q$ for all odd primes $q \neq 29$,
  and it is also annihilated by $\tau - 1$, where $\tau$ is induced
  by complex conjugation. We let $M'$ be the subgroup of~$M$
  annihilated by $\tau - 1$ and $T_q - \diamondop{q} - q$
  for $q = 3, 5, 7$. We find that
  \[ J_1(29)(\Q)_{\tors} \subseteq
     M' \cong \frac{\Z}{2\Z} \times \frac{\Z}{2\Z} \times \frac{\Z}{2\Z} \times
              \frac{\Z}{2^2\Z} \times \frac{\Z}{2^2\Z} \times
              \frac{\Z}{2^2 \cdot 3 \cdot 7 \cdot 43 \cdot 17837 \Z} .
  \]
  We can also compute the cuspidal subgroup~$C$ as the image in~$M$
  of the relative homology $H_1(X_1(29)(\C), \cusps, \Z)$ via
  its embedding into~$H_1(X_1(29)(\C), \Q)$. We obtain that
  \[ M' \subseteq C \cong \frac{\Z}{2^2\Z} \times \frac{\Z}{2^2\Z} \times
                          \frac{\Z}{2^2\Z} \times \frac{\Z}{2^2\Z} \times
                          \frac{\Z}{2^2 \cdot 3 \cdot 43 \cdot 17837 \Z} \times
                          \frac{\Z}{2^2 \cdot 3 \cdot 7^2 \cdot 43 \cdot 17837 \Z} .
  \]
  Finally, we have an explicit homomorphism $\Z[\cusps]^0 \to C$,
  where $\Z[\cusps]^0$ denotes the degree zero part of the free abelian
  group with basis the cusps of~$X_1(29)$. We know that the absolute
  Galois group of~$\Q$ fixes the $14$~cusps mapping to the cusp~$\infty$
  of~$X_0(29)$, whereas the remaining $14$~cusps are permuted cyclically
  via the action of the diamond operators. This allows us to determine
  \[ J_1(29)(\Q) = J_1(29)(\Q)_{\tors} = C^{\Gal_{\Q}}
                 \cong \frac{\Z}{2^2\Z} \times \frac{\Z}{2^2\Z} \times
                       \frac{\Z}{2^2 \cdot 3 \cdot 7 \cdot 43 \cdot 17837 \Z} ,
  \]
  and we can verify that this equals the subgroup generated by
  differences of rational cusps.
\end{proof}

\begin{remark}
  In our Magma code for the computations in the proof above, we rely
  only on linear algebra functionality (over $\Q$ and~$\Z$): we construct
  the relevant spaces ``by hand'' instead of using
  the built-in modular symbols functionality.
\end{remark}

Together with \cref{prop:rank0}, this immediately implies the following.

\begin{corollary} \label{cor:gen_by_rat_cusps}
  If $p$ is a prime such that $J_1(p)(\Q)$ is finite, then the
  latter group is generated by differences of rational cusps on~$X_1(p)$.
\end{corollary}

We will need the following result on the reduction mod~$2$.

\begin{proposition} \label{prop:beta_inj_small}
  If $p > 2$ is a prime such that $J_1(p)(\Q)$ is finite, then the reduction
  map $\red_2 \colon J_1(p)(\Q) = J_1(p)(\Q)_{\tors} \to J_1(p)(\F_2)$ is injective.
\end{proposition}

\begin{proof}
  Let $X$ be a curve over~$\Q$ with good reduction at~$2$,
  and let $J$ be its Jacobian variety.
  Then the kernel of the reduction map $J(\Q)_{\tors} \to J(\F_2)$
  is contained in the $2$-torsion subgroup~\cite{parent2}*{Lemme~1.7}.
  So the claim follows for all~$p$ such that $J_1(p)(\Q)$ is finite
  of odd order. For the remaining primes on our list,
  namely $p \in \set{17, 29, 31, 41}$, we check by an explicit computation
  that $J_1(p)(\Q)[2] \to J_1(p)(\F_2)$ is injective. This then implies
  the claim for these primes as well.

  We now describe this computation.
  By~\cref{cor:gen_by_rat_cusps}, we know that $J_1(p)(\Q)$ is generated
  by differences of rational cusps. The order of this group is known;
  see~\cite{CES}*{\S6.2.3 and Table~1} and note that the order for $p = 29$
  given there has to be divided by~$2^6$ to
  get the order of the group generated by differences of rational cusps;
  compare~\cref{thm:29cusp}.
  Sutherland~\cite{sutherland_table} provides equations for planar
  models of~$X_1(p)$ over~$\Q$ for the relevant values of~$p$.
  We use the reduction modulo~$2$ of this model
  to check that the subgroup of its Picard group generated by differences
  of its degree-$1$ places over~$\F_2$ (which correspond to the rational
  cusps under reduction mod~$2$) has the correct order. In fact, it suffices
  to check that the $2$-primary part of the group has the correct order.
  For $p = 17$, $29$, and~$31$, this only takes a few minutes;
  for $p = 41$ the computation of the Picard group of~$X_1(p)$ over~$\F_2$
  takes about eight hours (and $2.5$~gigabyte of memory).
\end{proof}

\begin{remark}
  If one does not want to wait for several hours for the computation
  for $p = 41$ to finish, one can alternatively use the intermediate
  curve~$X_H$ corresponding to $d = 4$ in the notation of~\cite{CES}
  (then $H$ has index~$4$).
  The predicted order of the $2$-primary part of~$J_H(\Q)$ equals that
  of~$J_1(p)(\Q)$. We check that the $2$-primary part of the subgroup
  of~$J_H(\F_2)$ generated by differences of the images of rational cusps has
  the correct size.
\end{remark}

\begin{remark}
  For $p = 17$, \cref{prop:beta_inj_small} together with the fact that
  the $\Q$-gonality of~$X_1(17)$ is~$4$ gives a simple alternative
  proof of the main result of~\cite{parent3} that $17 \notin S(3)$;
  see~\cref{cor:small_beta} below. (Note that $17 > (2^{3/2} + 1)^2$.)
\end{remark}

\begin{remark}
  The statement of~\cref{prop:beta_inj_small} is false for~$J_0(p)$.
  For example, $J_0(17)$ is the elliptic curve with Cremona label~$17a1$.
  It has $J_0(17)(\Q) \cong \Z/4\Z$, generated by the difference of the
  two cusps, but the reduction modulo~$2$ of a generator has only
  order~$2$.
\end{remark}

We now show that assumption~\eqref{rc_a} in~\cref{lem:crit_residue_classes}
is satisfied when $\ell = 2$ and $p$ is a rank zero prime such that the
$\Q$-gonality of~$X_1(p)$ is strictly larger than~$d$.
The \emph{$\Q$-gonality} of a curve~$X$ over~$\Q$ is the smallest degree
of a non-constant rational function on~$X$ defined over~$\Q$.

Recall the embedding $\iota \colon X_1(p)^{(d)} \to J_1(p)$ given by
fixing a base-point $c \in C_1(p)(\Q)$.

\begin{corollary} \label{cor:beta_inj_small}
  Let $d \ge 1$ be an integer.
  If $p > 2$ is a rank zero prime and the $\Q$-gonality of~$X_1(p)$
  is strictly larger than~$d$, then assumption~\eqref{rc_a}
  in~\cref{lem:crit_residue_classes} is satisfied for $\ell = 2$.
\end{corollary}

\begin{proof}
  The map $\iota \colon X_1(p)^{(d)}(\Q) \to J_1(p)(\Q)$ is injective
  when the $\Q$-gonality of~$X_1(p)$ exceeds~$d$, since otherwise
  there are two distinct $\Q$-rational effective divisors $D_1$ and~$D_2$
  of degree~$d$ that are linearly equivalent, which means that there
  is a rational function~$f$ on~$X_1(p)$ defined over~$\Q$ whose divisor
  is $D_1 - D_2$, hence $f$ has degree~$\le d$. This contradicts the
  condition on the $\Q$-gonality.
  \[ \xymatrix{ X_1(p)^{(d)}(\Q) \ar[d]^{\red_2} \ar[r]^-{\iota}
                  & J_1(p)(\Q) \ar[d]^{\red_2} \\
                X_1(p)^{(d)}(\F_2) \ar[r]^-{\iota}
                  & J_1(p)(\F_2)
              }
  \]
  By~\cref{prop:beta_inj_small}, the reduction map
  $\red_2 \colon J_1(p)(\Q) \to J_1(p)(\F_2)$ is injective as well,
  so $\red_2 \circ \iota = \iota \circ \red_2$ is injective, which
  implies that $\red_2$ is injective on~$X_1(p)^{(d)}(\Q)$.
\end{proof}

The following is an excerpt of~\cite{derickx_hoeij}*{Table~1}.
We write $\gon_{\Q}(X)$ for the $\Q$-gonality of a curve~$X$.
\[ \renewcommand{\arraystretch}{1.25}
  \begin{array}{|r|ccccccc|} \hline
                    p & 11 & 13 & 17 & 19 & 23 & 29 & 31 \\\hline
    \gon_{\Q}(X_1(p)) &  2 &  2 &  4 &  5 &  7 & 11 & 12 \\\hline
  \end{array}
\]

Also, it follows from~\cite{derickx_hoeij}*{Thm.~3} that
$\gon_{\Q}(X_1(p)) > 8$ for $p \in \{41, 47, 59, 71\}$.
We deduce the following.

\begin{corollary} \label{cor:small_beta}
  For $d$ and~$p$ as in the table below, assumption~\eqref{rc_a}
  in~\cref{lem:crit_residue_classes} is satisfied for $\ell = 2$.
  \begin{align*}
    d = 3 \colon & \quad p \in \set{17, 19, 23, 29, 31, 41, 47, 59, 71} \\
    d = 4 \colon & \quad p \in \set{19, 23, 29, 31, 41, 47, 59, 71} \\
    d = 5 \colon & \quad p \in \set{23, 29, 31, 41, 47, 59, 71} \\
    d = 6 \colon & \quad p \in \set{23, 29, 31, 41, 47, 59, 71} \\
    d = 7 \colon & \quad p \in \set{29, 31, 41, 47, 59, 71}
  \end{align*}
\end{corollary}

We now consider assumption~\eqref{rc_b} of~\cref{lem:crit_residue_classes}
for $p = 29, 31, 41$. We do this here rather than in~\cref{sec:ass_b},
since the computations we do to show that the assumption is satisfied
are closely related to those we do to establish~\cref{prop:beta_inj_small}.

\begin{lemma} \label{lem:29_31_41_alpha}
  For $p \in \set{29, 31, 41}$, $d \leq 7$ and $\ell = 2$,
  assumption~\eqref{rc_b} of~\cref{lem:crit_residue_classes}
  is satisfied.
\end{lemma}

\begin{proof}
  For $d \le 4$, we have that $p > (2^{d/2} + 1)^2$, and the claim follows
  from~\cref{lem:alpha_surj_simple}.
  For $(d,p) = (5,29)$, we observe that there is no elliptic curve
  over~$\F_{2^5}$ with $29$~points and that the cusps that are not
  images of rational cusps are not defined over~$\F_{2^5}$, so
  there are no points~$\bar{x}$ as in assumption~\eqref{rc_b}.

  In the other cases,
  \cref{cor:gen_by_rat_cusps} tells us that $J_1(p)(\Q)$ is generated
  by the differences of the rational cusps. This implies that the reduction
  mod~$2$ of any $\Q$-rational point of~$X_1(p)^{(d)}$ must map
  into the subgroup of~$J_1(p)(\F_2)$ that is generated by the differences
  of the images of the rational cusps. We verify
  that the points~$\bar{x}$ as in assumption~\eqref{rc_b}
  do not map into that subgroup, which by the above shows that these
  points are not in the image of the reduction map. This implies the claim.
  This computation is done together with the computations we do
  to prove~\cref{prop:beta_inj_small}.
\end{proof}

\begin{remark}
  In a similar way as in~\cref{prop:beta_inj_small}, we
  can use the following alternative approach for $p = 41$.
  There is no elliptic curve~$E$ over~$\F_{2^e}$ with $41 \mid \#E(\F_{2^e})$
  if $e \leq 7$ and $e \neq 5$.
  There is exactly one elliptic curve~$E$ over~$\F_{2^5}$
  with $\#E(\F_{2^5}) = 41$; this is the curve $y^2 + y = x^3 + x + 1$ already
  defined over~$\F_2$. Its automorphism group over~$\F_{2^5}$ is cyclic
  of order~$4$; we therefore obtain only $10 = (41-1)/4$ distinct
  $\F_{2^5}$-points on~$X_1(41)$ that are not cusps.
  Let $X_H$ be the intermediate curve between $X_1(41)$ and~$X_0(41)$
  with $H$ of index~$4$.
  Then $X_1(41) \to X_H$ is an \'etale cover of degree~$5$, and the ten
  $\F_{2^5}$-points on~$X_1(41)$ map to two $\F_2$-points on~$X_H$. In fact,
  $X_H(\F_2)$ consists of six points; four of them are cusps, and the other
  two are the ones just mentioned. It can be checked that these two points
  do not map into the subgroup generated by the differences of the four cusps,
  so that we can conclude in the same way as above.
\end{remark}


\section{Formal immersions} \label{sec:formal_imm}

When $p$ is not a rank zero prime, so that $J_1(p)(\Q)$ has positive rank,
then the reduction map $J_1(p)(\Q) \to J_1(p)(\F_\ell)$ is no longer
injective. This means that we need to find a more sophisticated argument
to verify assumption~\eqref{rc_a} of~\cref{lem:crit_residue_classes}.

As mentioned in the introduction,
one key idea here is to use a morphism $t \colon J_1(p) \to A$ of abelian
varieties over~$\Z_{(\ell)}$. We obtain the following commutative diagram.
\begin{equation} \label{eq:JA_diagram}
  \xymatrix{ X_1(p)^{(d)}(\Q) \ar[d]^{\red_\ell} \ar[r]^-{\iota}
              & J_1(p)(\Q) \ar[d]^{\red_\ell} \ar[r]^-{t}
              & A(\Q) \ar[d]^{\red_\ell} \\
             X_1(p)^{(d)}(\F_\ell) \ar[r]^-{\iota_\ell}
              & J_1(p)(\F_\ell) \ar[r]^-{t_\ell}
              & A(\F_\ell)
           }
\end{equation}
Let $\bar{x} \in X_1(p)^{(d)}(\F_\ell)$ be some point.
Assuming that $\red_\ell$ is injective on~$t(J_1(p)(\Q))$,
it will follow that the residue class of~$\bar{x}$ contains at most
one rational point, if we can show that the diagonal composition
$\red_\ell \circ t \circ \iota = t_\ell \circ \iota_\ell \circ \red_\ell$
is injective on the residue class of~$\bar{x}$.

The strategy for doing that is to take
$A$ such that $A(\Q)$ (or at least $t(J_1(p)(\Q))$) is finite; then the
reduction map on~$A(\Q)$ (or the image of~$t$) will be injective when
$\ell$ is odd; when $\ell = 2$, we can ensure that the reduction map
is injective on~$t(J_1(p)(\Q))$ by making sure that this image has odd
order. It then remains to show that $t \circ \iota$ is injective on the
residue class of any point $\bar{x} \in X_1(p)^{(d)}(\F_\ell)$.
To do this, we show that $t \circ \iota$ is a formal
immersion at each of the points~$\bar{x}$ as above.
We recall the definition below.

First, some notation. We write $\O_X$ for the structure sheaf of
a scheme~$X$, $\O_{X,x}$ for its local ring at a point~$x$ of~$X$,
and $\Ohat_{X,x}$ for the completion of the local ring with respect
to its maximal ideal~$\frm_{X,x}$.

\begin{definition} \label{def:formal_imm}
  Let $\phi \colon X \to Y$ be a morphism of noetherian schemes and let
  $x \in X$ be a point. Then $\phi$ is a \emph{formal immersion at $x$}
  if the induced local homomorphism on complete local rings
  \[ \hat{\phi}^* \colon \Ohat_{Y,\phi(x)} \to \Ohat_{X,x} \]
  is surjective.
\end{definition}

The relevant property of formal immersions for our purposes is the following;
this is (a consequence of)~\cite{parent1}*{Lemme~4.13}.

\begin{lemma} \label{lem:residue_class}
  Let $\phi \colon X \to Y$ be a morphism of noetherian schemes over~$\Z_{(\ell)}$
  that is a formal immersion at $x \in X(\F_\ell)$. Then $\phi$ induces
  an injective map on residue classes
  \[ \phi \colon \red_\ell^{-1}(x) \to \red_\ell^{-1}(\phi(x)) . \]
\end{lemma}

\begin{corollary} \label{cor:beta_inj}
  Let $d \in \Z_{\ge 1}$ and let $\ell \neq p$ be primes.
  Let $t \colon J_1(p) \to A$ be a morphism of abelian schemes
  over~$\Z_{(\ell)}$ such that
  \begin{enumerate}[\upshape(i)]
    \item \label{beta_inj1}
          $t(J_1(p)(\Q))$ is finite,
    \item \label{beta_inj2}
          $\ell > 2$ or $\# t(J_1(p)(\Q))$ is odd,
    \item \label{beta_inj3}
          $t \circ \iota$ is a formal immersion at all
          $\bar{x} \in X_1(p)^{(d)}(\F_\ell)$
          that are sums of images of rational cusps on~$X_1(p)$.
  \end{enumerate}
  Then assumption~\eqref{rc_a} of~\cref{lem:crit_residue_classes} is
  satisfied.
\end{corollary}

\begin{proof}
  Note that $X_1(p)$ and~$J_1(p)$ have good reduction at~$\ell$, hence
  $J_1(p)$ can be considered as an abelian scheme over~$\Z_{(\ell)}$.

  Let $x, x' \in X_1(p)^{(d)}(\Q)$ be in the residue class of a
  point~$\bar{x}$ that is a sum of images of rational cusps and
  write $y = t(\iota(x))$, $y' = t(\iota(x'))$. Since $x$ and~$x'$
  are in the same residue class, the same is true of $y$ and~$y'$.
  It follows from conditions \eqref{beta_inj1} and~\eqref{beta_inj2}
  that $\red_\ell \colon t(J_1(p)(\Q)) \to A(\F_\ell)$ is injective,
  which implies that $y = y'$.
  By condition~\eqref{beta_inj3} and~\cref{lem:residue_class},
  $t \circ \iota$ is injective on the residue class of~$\bar{x}$,
  which finally shows that $x = x'$.
\end{proof}

\begin{remark}
  If $\ell = 2$ and we take commuting $t_1, t_2 \in \End_{\Q}(J_1(p))$
  such that $t_1(J_1(p)(\Q))$ is finite and $t_2$ kills the $2$-torsion
  subgroup of~$J_1(p)(\Q)$, then the conclusion of~\cref{cor:beta_inj}
  holds for $t = t_1 t_2$ also when $\#t(J_1(p)(\Q))$ is even (assuming
  condition~\eqref{beta_inj3} is satisfied); see ~\cite{parent2}*{Thm.~1.10}.
  Writing $A_1 = \im(t_1)$ and taking $A = \im(t)$ without loss of
  generality, we have the following commutative diagram.
  \[ \xymatrix{ X_1(p)^{(d)}(\Q) \ar[d]^{\red_2} \ar[r]^-{\iota}
                 & J_1(p)(\Q) \ar[d]^{\red_2} \ar[r]^-{t_1}
                 & A_1(\Q) \ar[d]^{\red_2} \ar[r]^-{t_2}
                 & A(\Q) \ar[d]^{\red_2} \\
                X_1(p)^{(d)}(\F_2) \ar[r]^-{\iota}
                 & J_1(p)(\F_2) \ar[r]^-{t_1}
                 & A_1(\F_2) \ar[r]^-{t_2}
                 & A(\F_2)
              }
  \]
  Take $x, x' \in X_1(p)^{(d)}(\Q)$ with the same reduction~$\bar{x}$ mod~$2$,
  such that $\bar{x}$ is a sum of images of rational cusps.
  Then $t_1(\iota(x') - \iota(x))$ is in the kernel of reduction mod~$2$
  of~$A_1(\Q)$, which (since $A_1(\Q)$ is finite) consists of $2$-torsion
  points, so $t(\iota(x')) = t(\iota(x))$ by the assumption on~$t_2$.
  We can then conclude as in the proof above.
\end{remark}

In our intended application, the set $X_1(p)^{(d)}(\F_\ell)$ can be
quite large: the curve $X_1(p)$ has $(p-1)/2$ $\Q$-rational cusps;
assuming that they account for all of~$X_1(p)^{(d)}(\F_\ell)$, the latter
set has $\binom{(p-1)/2 + d - 1}{d}$ elements. \cref{cor:beta_inj}
requires us to check that $t \circ \iota$ is a formal immersion at
each of these points. To reduce the necessary computational effort,
we now show how we can use curves intermediate between $X_1(p)$
and~$X_0(p)$ that have fewer cusps.

\begin{corollary} \label{cor:intermediate_curve}
  Let $d \in \Z_{\ge 1}$ and let $\ell \neq p$ be primes.
  Let $X_H$ be an intermediate curve between $X_1(p)$ and~$X_0(p)$.
  Fix $x_0 \in X_H^{(d)}(\Q)$ and define $\iota_H \colon X_H^{(d)} \to J_H$
  using~$x_0$ as base-point. Let $t \colon J_H \to A$ be a morphism of
  abelian schemes over~$\Z_{(\ell)}$ such that
  \begin{enumerate}[\upshape(i)]
    \item \label{int_curve1}
          $t(J_H(\Q))$ is finite,
    \item \label{int_curve2}
          $\ell > 2$ or $\# t(J_H(\Q))$ is odd,
    \item \label{int_curve3}
          $t \circ \iota_H$ is a formal immersion at all
          $\bar{x}_H \in X_H^{(d)}(\F_\ell)$
          that are sums of images of rational cusps on~$X_1(p)$.
  \end{enumerate}
  Then assumption~\eqref{rc_a} of~\cref{lem:crit_residue_classes} is
  satisfied.
\end{corollary}

\begin{proof}
  Let $\bar{x} \in X_1(p)^{(d)}(\F_\ell)$ be a sum of images of rational
  cusps and take two points $x, x' \in X_1(p)^{(d)}(\Q)$ in the residue
  class of~$\bar{x}$, where we take $x$ to be the unique sum of rational
  cusps in this residue class.
  Write $x_H, x'_H$ for their images in~$X_H^{(d)}(\Q)$.
  Then $\bar{x}_H := \red_\ell(x_H) = \red_\ell(x'_H)$ is a sum of images
  of rational cusps on~$X_1(p)$. Arguing as in the proof
  of~\cref{cor:beta_inj}, we see that $x'_H = x_H$; in particular, $x'_H$
  is a sum of images of rational cusps on~$X_1(p)$, since this is true
  for~$x_H$. The set of rational cusps on~$X_1(p)$ is the full preimage
  of the cusp $\infty \in X_0(p)(\Q)$. This implies that all points
  in~$X_1(p)^{(d)}(\Q)$ that are preimages of~$x_H$ under the obvious
  map are sums of rational cusps. So $x'$ is a sum of rational cusps
  as well. But $\red_\ell$ is injective on sums of rational cusps
  (since reduction mod~$\ell$ is injective on cusps;
  see~\cite{deligne_rapoport}*{Thm.~IV.3.4}), hence $x' = x$.
\end{proof}


\section{Computational verification of assumption \eqref{rc_a}}
\label{sec:ass_a}

We use~\cref{cor:intermediate_curve} to show that
assumption~\eqref{rc_a} of~\cref{lem:crit_residue_classes} holds for
the relevant pairs~$(d,p)$.
To verify the assumptions of~\cref{cor:intermediate_curve}, we need to do
essentially two things: we have to find a suitable morphism~$t$ of abelian
schemes that satisfies conditions \eqref{int_curve1} and~\eqref{int_curve2},
and we have to check that $t \circ \iota$ is a formal immersion at all
points in~$\bar{x}_H \in X_H^{(d)}(\F_\ell)$ that are sums of images
of rational cusps on~$X_1(p)$.

To satisfy condition~\eqref{int_curve1}, we take a morphism~$t$
that factors through the winding quotient~$J^{\bfe}_H$; then $t(J_H(\Q))$
is contained in the image of $J^{\bfe}_H(\Q)$ under a morphism of abelian
varieties. Since by~\cref{thm:winding} $J^{\bfe}_H(\Q)$ is finite,
$t(J_H(\Q))$ is finite as well. One possibility is to take the projection
$J_H \to J^{\bfe}_H$. If we choose $\ell \ge 3$, then
condition~\eqref{int_curve2} is also satisfied.
This was used for~$J_0(p)$ with $p$~prime and an~$\ell$
that depends on~$p$ in the argument of~\cite{merel},
and is used for~$J_0(p^n)$ with $\ell = 3$ or~$5$ in the argument
of~\cite{parent1}.
The proof of Oesterl\'e's bound uses $\ell = 3$; see~\cref{sec:oesterle}.

If we take for~$t$ an element
of the Hecke algebra~$\T \subseteq \End_{\Q}(J_H)$,
then the condition for $t$ to factor via the winding quotient is that
$t \cdot \Ann(\bfe) = 0$ in $\T$. We obtain such~$t$ as follows.
This is essentially~\cite{parent2}*{Lemme~1.9}; we extend the statement
slightly by removing the condition that the characteristic polynomial
of~$t_0$ (acting on the space of cusp forms) is squarefree.

\begin{proposition} \label{prop:t1}
  Let $t_0 \in \T$ with factored characteristic
  polynomial $P(X) = \prod_{i=1}^n P_i(X)^{e_i}$
  with respect to its action on~$H^0(X_H, \Omega^1)$. Set
  \[ I := \set{i \in \set{1,\ldots,n} \mid (P/P_i)(t_0) \cdot \bfe = 0
                                           \text{\ or\ } e_i \ge 2} ;
  \]
  then $t_1(t_0) := \prod_{i \in I} P_i^{e_i}(t_0)$ is such that
  $t_1(t_0) \cdot \Ann(\bfe) = 0$.
\end{proposition}

\begin{proof}
  The proof is basically the same as that in~\cite{parent2}*{\S2.5},
  noting that the factors $P_i^{e_i}(t_0)$ with $e_i \ge 2$ in the
  product defining~$t_1(t_0)$ are used to kill any factor of the
  Hecke algebra for which we cannot simply decide whether it is
  contained in~$\Ann(\bfe)$.
\end{proof}

We note that we can compute~$P(X)$ and test the condition
$(P/P_i)(t_0) \cdot \bfe = 0$ explicitly using modular symbols, so we can
determine~$t_1(t_0)$ explicitly for any given~$t_0$. We see that~$t_1(t_0)$
satisfies condition~\eqref{int_curve1} for every $t_0 \in \T$.

To satisfy condition~\eqref{int_curve2} when $\ell = 2$, we
use~\cref{prop:ann_rat_tors}, which implies that for $q$ an odd prime
not dividing~$N$,
$T_q - \diamondop{q} - q$ kills the rational torsion subgroup of~$J_H$.
Combining this with~\cref{prop:t1}
gives the following version of ``Parent's trick''.

\begin{corollary} \label{cor:trick}
  Let $X_H$ be an intermediate curve between $X_1(p)$ and~$X_0(p)$.
  Let $t_0 \in \T$ and let $q \neq p$ be an odd prime. Then
  \[ t := t_1(t_0) \cdot (T_q - \diamondop{q} - q) \in \T , \]
  considered as an element of~$\End_{\Q}(J_H)$,
  satisfies conditions \eqref{int_curve1} and~\eqref{int_curve2}
  of~\cref{cor:intermediate_curve} for $\ell = 2$.
  If $X_H = X_0(p)$ and $p \not\equiv 1 \bmod 8$, then $t := t_1(t_0)$
  satisfies both conditions.
\end{corollary}

\begin{proof}
  By~\cref{prop:t1} and the discussion preceding it, $t_1(t_0)$ satisfies
  condition~\eqref{int_curve1}. Obviously this condition still holds
  after composing $t_1(t_0)$ with some further morphism.
  By~\cref{prop:ann_rat_tors}, the factor $T_q - \diamondop{q} - q$
  kills the torsion
  in $t_1(t_0)(J_H(\Q)) \subseteq J_H(\Q)_{\tors}$, which implies that
  $t(J_H(\Q)) = 0$, so that condition~\eqref{int_curve2} also holds.

  It is known that $J_0(p)(\Q)_{\tors}$ is cyclic of order
  $(p-1)/\gcd(p-1, 12)$, generated by the difference of the two
  (rational) cusps; see~\cite{mazur1}*{Thm.~1}. This implies that
  the rational torsion group of~$J_0(p)$ has odd order when
  $p \not\equiv 1 \bmod 8$, and so condition~\eqref{int_curve2}
  is automatically satisfied.
\end{proof}

We still need a way of verifying condition~\eqref{int_curve3}
of~\cref{cor:intermediate_curve}. This is provided by the following
version of ``Kamienny's criterion'' as given
in~\cite{parent2}*{Thm.~1.10, Prop.~2.7}.
Parent states this criterion for~$X_1(p)$ in place of~$X_H$, but the
generalization is immediate.

\begin{proposition} \label{prop:kamiennys_criterion}
  Let $H \subseteq (\Z/p\Z)^\times/\set{\pm 1}$ be a subgroup.
  Let $\ell \neq p$ be a prime and consider $t = t_1(t_0)$ as
  in~\cref{prop:t1} when $\ell$ is odd, or $t$ as in~\cref{cor:trick}
  when $\ell = 2$. Then $t \circ \iota$ is a formal immersion at all
  $\bar{x}_H \in X_H^{(d)}(\F_\ell)$
  that are sums of images of rational cusps on~$X_1(p)$,
  if for all partitions $d = n_1 + \ldots + n_m$ with $n_1 \ge \dots \ge n_m$
  and all $m$-tuples
  $(d_1 = 1, d_2, \ldots, d_m)$ of integers representing pairwise distinct
  elements of~$H$, the $d$ Hecke operators
  \begin{equation} \label{eq:kamienny}
    \bigl(T_i \diamondop{d_j} t\bigr)_{\substack{j = 1, \ldots, m \\ i = 1, \ldots, n_j}}
  \end{equation}
  are $\F_\ell$-linearly independent in~$\T \otimes \F_\ell$, where
  $\T$ is considered as a subalgebra of~$\End_{\Q}(J_H)$.
\end{proposition}

We note that we can check the criterion for any given~$t$ by a computation
with modular symbols.

This criterion was first established by Kamienny in~\cite{kamienny1}
for~$X_0(p)$. In this case, the condition simplifies to

\emph{the $d$ Hecke operators $T_1 t, T_2 t, \ldots, T_d t$
are $\F_\ell$-linearly independent in $\T \otimes \F_\ell$.}

Implementing the criterion implied by~\cref{cor:trick,prop:kamiennys_criterion}
for~$X_0(p)$ and running the resulting code gives the following.
We take as base-point for~$\iota$ the point given by $d$~times
the cusp~$\infty$ on~$X_0(p)$. Note that $2281 = \floor{(1 + 3^{7/2})^2}$;
larger primes will be dealt with using Oesterl\'e's bound.

\begin{lemma} \label{lem:formal_imm_X0}
  For each of the following choices of $3 \le d \le 7$ and a prime~$p$,
  there is $t \in \End_{\Q}(J_0(p))$ as in~\cref{cor:trick} for $\ell = 2$
  such that
  $t \circ \iota \colon X_0(p)^{(d)}_{\Z_{(2)}} \to J_0(p)_{\Z_{(2)}}$
  is a formal immersion at the point of~$X_0(p)^{(d)}(\F_2)$ corresponding
  to $d$~times the cusp~$\infty$:
  \begin{align*}
    d = 3 \colon & \quad 47 \le p \le 2281, \quad p \neq 73, 79; \\
    d = 4 \colon & \quad p \in \set{47, 59, 71, 83, 89}
                         \quad\text{or}\quad 103 \le p \le 2281; \\
    d = 5 \colon & \quad p \in \set{59, 71, 83}
                         \quad\text{or}\quad 103 \le p \le 2281; \\
    d = 6 \colon & \quad p \in \set{71, 107}
                         \quad\text{or}\quad 127 \le p \le 2281, \; p \neq 193; \\
    d = 7 \colon & \quad p = 131 \quad\text{or}\quad
                         139 \le p \le 2281, \; p \neq 157, 193.
  \end{align*}
\end{lemma}

\begin{proof}
  We try $t_0 = T_n$ for
  $2 \le n \le 60$, and when $p \equiv 1 \bmod 8$, we try for each~$t_0$
  the additional factor $T_q - (q+1)$ for primes $3 \le q \le 20$
  until either the criterion is satisfied or else all combinations
  are exhausted. (Actually, $n \le 14$ and $q \in \set{3,5}$ would be
  enough, as the computation reveals). The computation took about $1.5$~hours.
  We note that to exclude $p = 163$ for $d = 7$,
  we actually needed the statement of~\cref{cor:trick} that \hbox{$t = t_1(t_0)$}
  is sufficient when $p \not\equiv 1 \bmod 8$
  (which also helps to speed up the computation, since it eliminates
  the inner loop over~$q$).
  For $p = 431$ and~$d = 7$, taking $t_0 = T_n$ does not seem to work.
  We tried random linear combinations of the first few Hecke operators
  and were successful with $t_0 = T_2 + T_3 - T_7$.
\end{proof}

For the remaining primes~$p$ of interest for any given degree~$d$,
we use the criterion on an intermediate curve~$X_H$; we try the various
groups~$H$ ordered by increasing index in~$(\Z/p\Z)^\times/\set{\pm 1}$,
since smaller index means that we have to deal with smaller objects,
leading to a faster computation.

If we were to use the criterion of~\cref{prop:kamiennys_criterion}
literally, then we would have to run through a potentially very large
number of partitions of~$d$ combined with choices of~$d_j$. We use
the following trick to speed up the computation.

\begin{lemma} \label{lemma:faster_criterion}
  Let $H \subseteq (\Z/p\Z)^\times/\set{\pm 1}$ be a subgroup.
  Let $\ell \neq p$ be a prime, $d$ an integer and $t \in \T$,
  viewed as an endomorphism of~$J_H$.
  Let $D \subseteq \Z$ be a set of representatives of the cosets of~$H$
  with $1 \in D$. Define
  the set
  \[ I := \set{(1,i) \mid 1 \le i \leq d}
            \cup \set{(k,i) \mid 1 \leq i \leq \floor{d/2}, 1 \neq k \in D}.
  \]
  Suppose that
  there is no $\F_\ell$-linear dependence among at most~$d$ of the
  images in $\T \otimes \F_\ell$ of the elements
  $t_{(k,i)} := T_i \diamondop{k} t$ for $(k,i) \in I$, where
  we consider $\T$ as a subalgebra of~$\End_{\Q}(J_H)$.
  Then the criterion of~\cref{prop:kamiennys_criterion} is satisfied.
\end{lemma}

\begin{proof}
  Assume the criterion fails. Then there is a partition
  $d = n_1 + \dots + n_m$ with $n_1 \ge \dots \ge n_m$
  and there are $d_1 = 1$, $d_2, \ldots, d_m \in D$ pairwise distinct
  such that the $d$~operators
  $T_i \diamondop{d_j} t$ for $1 \le j \le m$, $1 \le i \le n_j$
  are linearly dependent in~$\T \otimes \F_\ell$. But these operators
  are all of the form~$t_{(k,i)}$ (note that $n_j \le \floor{d/2}$
  for $j \ge 2$), so this would produce a linear
  dependence mod~$\ell$ among $d$ of the~$t_{(k,i)}$; this is a
  contradiction.
\end{proof}

When implementing this, we can in addition look at each linear
relation of weight at most~$d$ between the elements in the lemma
and check if it is indeed of the ``forbidden'' form as given
in~\cref{prop:kamiennys_criterion}. In the cases
of interest, the relation space has low enough dimension to allow
for the enumeration of all relations and performing this check.
We use algorithms for binary linear codes that are included
in Magma to do this efficiently.

We obtain the following result.

\begin{lemma} \label{lem:formal_imm}
  For each of the following choices of $3 \le d \le 7$ and a prime~$p$,
  there is a subgroup~$H$ of~$(\Z/p\Z)^\times/\set{\pm 1}$ and
  $t \in \End_{\Q}(J_H)$ as in~\cref{cor:trick} for $\ell = 2$
  such that
  $t \circ \iota \colon X^{(d)}_{H,\Z_{(2)}} \to J_{H,\Z_{(2)}}$
  is a formal immersion at all points of~$X_H^{(d)}(\F_2)$ that are
  sums of images of rational cusps on~$X_1(p)$:
  \begin{align*}
    d = 3 \colon & \quad 19 \le p \le 2281; \\
    d = 4 \colon & \quad 19 \le p \le 2281, \quad p \neq 29; \\
    d = 5 \colon & \quad 23 \le p \le 2281, \quad p \neq 29; \\
    d = 6 \colon & \quad 23 \le p \le 2281, \quad p \neq 29; \\
    d = 7 \colon & \quad 37 \le p \le 2281.
  \end{align*}
\end{lemma}

\begin{proof}
  For each pair~$(d,p)$
  that is not covered by~\cref{lem:formal_imm_X0}, we check the
  criterion of~\cref{lemma:faster_criterion} for subgroups~$H$
  by increasing index. For each~$H$, we again try $t_0 = T_n$
  for $2 \le n \le 60$ and the second factor given by primes
  $3 \le q \le 20$. The most involved computation is for $d = 7$
  and $p = 107$, where we have to take the trivial subgroup~$H$
  corresponding to~$J_1(107)$; this computation took about
  $35$~minutes. Most of the other cases just take a few seconds,
  a small number of them a few minutes.
\end{proof}

\cref{prop:ass_a} now follows from~\cref{lem:formal_imm,cor:small_beta}.


\section{A proof of Oesterl\'e's bound} \label{sec:oesterle}

The purpose of this section is to provide a proof of Oesterl\'e's
bound~\eqref{eq:oesterle} and thus close a gap in the literature.
Oesterl\'e gives a proof in his notes~\cite{oesterle}, which have
been available to the people working in the field, but a proof has
never appeared in print. The proof below is based on these notes,
which Oesterl\'e kindly provided to us; in particular, we do not
claim originality for anything in this section: the ideas are all
Oesterl\'e's. We will use results that are available in the literature
by now to simplify the exposition in some places. We state the result
of this section as a theorem.

\begin{theorem}[Oesterl\'e] \label{thm:oesterle}
  Let $d \ge 3$. If $p > (3^{d/2} + 1)^2$ is a prime, then $p \notin S(d)$.
\end{theorem}

We can restrict to $d \ge 3$ here, since the cases $d = 1$ and~$d = 2$
have been dealt with by Mazur and Kamienny, respectively.

We will work with $\ell = 3$. By \cref{lem:alpha_surj_simple},
assumption~\eqref{rc_b} of~\cref{lem:crit_residue_classes} is always
satisfied when $p > (3^{d/2} + 1)^2$. So it is sufficient to show
that assumption~\eqref{rc_a} of~\cref{lem:crit_residue_classes} holds.
This in turn is done by using the formal immersion criterion via
the winding quotient of~$J_0(p)$. For sufficiently large~$d$, this
follows from the following result.

\begin{proposition} \label{prop:large_d}
  If $d \ge 3$ and $p \ge 65 (2 d)^6$ is a prime, then the map
  \[ f_{d,p} \colon X_0(p)^{(d)} \stackrel{\iota}{\to} J_0(p) \to J_0^{\bfe}(p) \]
  is a formal immersion at the point $\bar{x} \in X_0(p)^{(d)}(\F_3)$ that is
  the reduction mod~$3$ of $d$~times the cusp~$\infty$ on~$X_0(p)$.

  In particular, \cref{thm:oesterle} holds for $d \ge 26$.
\end{proposition}

\begin{proof}
  The first statement is a consequence of~\cite{parent1}*{Thm.~1.8 and Prop.~1.9}.
  Since $65 (2d)^6 < (3^{d/2} + 1)^2$ when $d \ge 26$, the statement
  of~\cref{thm:oesterle} follows for such~$d$ by the discussion above.
\end{proof}

Oesterl\'e proves a similar statement with a slightly worse bound on~$p$;
Parent uses the same underlying approach.

In principle, \cref{prop:large_d} reduces the proof of~\cref{thm:oesterle}
to a finite problem: for each $3 \le d \le 25$ and each prime~$p$ such that
$(3^{d/2} + 1)^2 < p < 65 (2d)^6$, we have to check that the map
in~\cref{prop:large_d} is a formal immersion at the relevant point,
which can be done via Kamienny's criterion given
in~\cref{prop:kamiennys_criterion}. However, the primes we would have
to deal with in this way get much too large and there are way too many
of them to make this practical. So instead, we need a criterion that
allows us to deal with all (or many) of these primes at the same time.

One idea that Oesterl\'e uses here (and also to prove a statement similar
to~\cref{prop:large_d} above) is to make use of the intersection pairing
on~$H_1(X_0(p)(\C), \Z)$, which is an alternating perfect pairing into~$\Z$.
We will denote this pairing by~$\bullet$.

We will use the following version of Kamienny's criterion.
Recall the winding element~$\bfe \in H_1(X_0(p)(\C), \Q)$
from~\cref{def:winding}.

\begin{proposition} \label{prop:kamienny0}
  The map $f_{d,p}$ as in~\cref{prop:large_d} is a formal immersion
  at~$\bar{x}$ if (and only if) the images of $T_1 \bfe, \ldots, T_d \bfe$
  in~$\T\bfe/3\T\bfe$ are linearly independent over~$\F_3$.
\end{proposition}

\begin{proof}
  This is~\cite{parent1}*{Thm.~4.18} for $l = 3$.
\end{proof}

To use the intersection pairing, we have to move the elements
$T_n \bfe \in H_1(X_0(p)(\C), \Q)$ into~$H_1(X_0(p)(\C), \Z)$.
The Hecke operator $T_2 - 3$ sends~$\bfe$ into~$H_1(X_0(p)(\C), \Z)$,
since the action of~$T_2$, viewed as a correspondence on~$X_0(p)$,
multiplies the cusps $0$ and~$\infty$ by~$3$, so that the boundary
of $-(T_2 - 3) \cdot \set{0, \infty}$ is zero.
In the same way, we see that
$(T_n - \sigma_1(n)) \bfe \in H_1(X_0(p)(\C), \Z)$
when $n < p$; here $\sigma_1(n)$ denotes the sum of (positive)
divisors of~$n$. (This is true in general when $p \nmid n$;
when $p \mid n$, one has to replace $\sigma_1(n)$ with the sum
of divisors not divisible by~$p$.)

\begin{corollary} \label{cor:kamienny_Z}
  If $p > (3^{d/2} + 1)^2$ and the images of
  \[ (T_2 - 3) T_1 \bfe, \ldots, (T_2 - 3) T_d \bfe \]
  in~$H_1(X_0(p)(\C), \F_3)$ are linearly independent over~$\F_3$,
  then $p \notin S(d)$.
\end{corollary}

\begin{proof}
  We show that $f_{d,p}$ is a formal immersion at~$\bar{x}$, which
  implies the claim. Assume that this is not the case.
  By~\cref{prop:kamienny0}, there are integers
  $\lambda_1, \ldots, \lambda_d$, not all divisible by~$3$, such that
  $\lambda_1 T_1 \bfe + \dots + \lambda_d T_d \bfe \in 3 \T \bfe$.
  Multiplying by~$T_2 - 3$, this gives
  \[ \lambda_1 (T_2 - 3) T_1 \bfe + \dots + \lambda_d (T_2 - 3) T_d \bfe
       \in 3 (T_2 - 3) \T \bfe \subset 3 H_1(X_0(p)(\C), \Z) \,,
  \]
  with all terms on the left contained in~$H_1(X_0(p)(\C), \Z)$.
  Reducing this relation mod~$3$ shows that the images of
  $(T_2 - 3) T_1 \bfe, \ldots, (T_2 - 3) T_d \bfe$
  in~$H_1(X_0(p)(\C), \F_3)$ are linearly dependent.
\end{proof}

We now define the following Hecke operators.

\begin{definition}
  Let $n \ge 1$. We set
  \[ T'_n = \sum_{m \mid n} \mu\bigl(\tfrac{n}{m}\bigr) T_m \,, \]
  where $\mu$ is the M\"obius function, and
  \[ L_n = T'_{2n} - 2 T'_n \,. \]
\end{definition}

Then $T_n - \sigma_1(n) = \sum_{m \mid n} (T'_m - m)$.
Using the relations
\[ T_2 T_m = \begin{cases}
               \hfill T_{2m}, & \text{if $m$ is odd,} \\
               T_{2m} + 2 T_{m/2}, & \text{if $m$ is even,}
             \end{cases}
\]
we find that
\[ (T_2 - 3) T'_n = \begin{cases}
                      \hfill L_n, & \text{if $n$ is odd,} \\
                      L_n - L_{n/2}, & \text{if $n$ is even.}
                    \end{cases}
\]

\begin{corollary} \label{cor:kamienny_L}
  If $p > (3^{d/2} + 1)^2$ and the images of
  \[ L_1 \bfe, \ldots, L_d \bfe \]
  in~$H_1(X_0(p)(\C), \F_3)$ are linearly independent over~$\F_3$,
  then $p \notin S(d)$.
\end{corollary}

\begin{proof}
  The relations deduced above show that the $\Z$-submodule of~$\T$
  generated by $L_1, \ldots, L_d$ is the same as the $\Z$-submodule
  generated by $(T_2 - 3) T_1, \ldots, (T_2 - 3) T_d$. Now
  use~\cref{cor:kamienny_Z}.
\end{proof}

We now introduce notation for certain modular symbols,
following~\cite{merel}*{Section~2}.
If $\gamma = \smm{a}{b}{c}{d} \in \SL_2(\Z)$,
then the modular symbol $\set{\gamma 0, \gamma \infty}$ depends only
on the coset $\Gamma_0(p) \gamma$, which in turn depends only on
the image of $\tfrac{c}{d}$ in $\P^1(\F_p)$. We denote this modular
symbol by~$\xi(\tfrac{c}{d})$. If $k$ is an integer coprime with~$p$,
then $\xi(k) = \set{0, \tfrac{1}{k}} \in H_1(X_0(p)(\C), \Z)$,
since the cusp~$\tfrac{1}{k}$ is $\Gamma_0(p)$-equivalent to~$0$.

The following result is crucial; we defer its proof
until later and first show how~\cref{thm:oesterle}
can be deduced from it with some computation.
For $M \ge 3$ an odd integer, we define
\[ \eps_M \colon (\Z/M\Z)^\times \to \set{0,1} \]
so that $\eps_M(a + M\Z) = 0$ if $1 \le a < M/2$ and
$\eps_M(a + M\Z) = 1$ if $M/2 < a < M$. We extend $\eps_M$
to a map on all rational numbers~$a/b$ with numerator and denominator
coprime to~$M$ by applying it to the image of~$a/b$ in~$(\Z/M\Z)^\times$.

\begin{lemma} \label{lem:inters}
  Let $d \ge 1$ be an integer, let $M \ge 3$ be an odd integer
  and let $p > 2dM$ be a prime.
  Let $u \in \Z$ be such that $p u \equiv 1 \bmod M$.
  Then for $a$ coprime to~$M$ and $1 \le n \le d$, we have
  \[ L_n \bfe \bullet \set{0, \frac{a}{M}} = \eps_M(na) - \eps_M(nu/a) . \]
\end{lemma}

\begin{corollary}[\cite{oesterle}*{Prop.~8}] \label{cor:oesterle_matrix}
  Let $d$ and~$M$ be as in~\cref{lem:inters} and fix $u \in \Z$
  coprime with~$M$. If the matrix
  \[ \bigl(\eps_M(na) - \eps_M(nu/a)\bigr)_{1 \le n \le d, a \in (\Z/M\Z)^\times} , \]
  with entries taken in~$\F_3$, has rank~$d$, then $p \notin S(d)$
  for all primes
  \[ p > \max\{2dM, (3^{d/2} + 1)^2\}
     \quad\text{such that $p u \equiv 1 \bmod M$.}
  \]
\end{corollary}

\begin{proof}
  By~\cref{lem:inters}, the matrix entries are the intersection numbers,
  taken mod~$3$, between $L_n \bfe$ and $\set{0,\tfrac{a}{M}}$
  in $H_1(X_0(p)(\C), \Z)$, when $p$ is a prime as in the statement.
  So when the matrix has rank~$d$, this implies that
  $L_1 \bfe, \ldots, L_d \bfe$ are linearly independent mod~$3$,
  and the claim follows from~\cref{cor:kamienny_L}.
\end{proof}

\begin{proof}[Proof of~\cref{thm:oesterle}]
  The following table gives, for each $3 \le d \le 25$, a value of~$M$
  as in~\cref{cor:oesterle_matrix} such that the matrix above has rank~$d$
  for all $u \in (\Z/M\Z)^\times$. By~\cref{cor:oesterle_matrix}, this
  proves~\cref{thm:oesterle} for all $p > \max\{2dM, (3^{d/2}+1)^2\}$.
  \[ \renewcommand{\arraystretch}{1.2}
     \begin{array}{|c||*{12}{c|}} \hline
       d &  3 &  4 &  5 &  6 &  7 &  8 &  9 & 10 & 11 & 12 & 13 & 14 \\\hline
       M & 29 & 37 & 41 & 43 & 47 & 47 & 53 & 53 & 53 & 61 & 73 & 73 \\\hline
       \multicolumn{13}{c}{\strut} \\[-6pt]\hline
       d & 15 & 16 & 17 & 18 & 19 & 20 & 21 & 22 & 23 & 24 & 25 &    \\\hline
       M & 79 & 79 & 89 & 89 & 89 & 101 & 101 & 109 & 109 & 109 & 127 & \\\hline
     \end{array}
  \]
  Note that $2dM < (3^{d/2} + 1)^2$ for $d \ge 6$. We have already
  verified the formal immersion criterion (with $\ell = 2$) for the primes
  between $(3^{d/2} + 1)^2$ and~$2dM$ for $3 \le d \le 5$
  in~\cref{lem:formal_imm}, which implies $p \notin S(d)$ for these
  primes as well.
\end{proof}

\begin{remark}
  Oesterl\'e deals with the remaining primes~$p$ for $3 \le d \le 5$
  by computing the intersection products $I_n \bfe \bullet \xi(k)$
  for $1 \le k \le p-1$ and $1 \le n \le d$, where
  \hbox{$I_1 = (p-1)/\gcd(p-1,12)$} is the order of~$J_0(p)(\Q)_{\tors}$
  and $I_n = T'_n - n$ for $n \ge 2$, and verifying that the resulting
  matrix has rank~$d$ (even when reduced modulo any prime $\ell \ge 3$).
  This works for all cases except $p = 43$ and $p = 73$ for $d = 3$.
  For $p = 73$, he has a separate argument, whereas he does not mention
  $p = 43$ further, even though the maximal~$d$ for which the rank condition
  is satisfied is $d = 2$ according to the table at the end
  of~\cite{oesterle}*{Section~7}.
\end{remark}

From now on until the end of this section, the degree~$d$ of the field
of definition of the elliptic curves will be irrelevant. We will therefore
feel free to use ``$d$'' as a local variable as in the definition
of~$M_n$ below, and hope that this will not lead to confusion.

It remains to prove~\cref{lem:inters}. We follow Oesterl\'e's note
quite closely here (modulo some changes of notation). We remark that
\cref{cor:inters_large_p,cor:inters_aM} are in a separate file that
Oesterl\'e made available to Bas Edixhoven and the first author of this paper.

We begin with a result that expresses
$(T_n - \sigma_1(n)) \bfe$ for $n < p$ in terms of modular symbols.

\begin{lemma}[\cite{oesterle}*{Cor.~2 of Prop.~10}] \label{lem:Te_modsym}
  For $n < p$, we have in~$H_1(X_0(p)(\C), \Z)$ that
  \[ (T_n - \sigma_1(n)) \bfe = -\sum_{(a,b,c,d) \in M_n} \xi(\tfrac{c}{d}) , \]
  where
  \[ M_n = \{(a,b,c,d) \in \Z^4 : a > b \ge 0, \; d > c > 0, \; ad - bc = n\} . \]
\end{lemma}

\begin{proof}
  This is~\cite{merel}*{Lemme~2}, taking into account that both sides
  are contained in~$H_1(X_0(p)(\C), \Z)$.
\end{proof}

We need a formula for the intersection product. We define the following
function on~$\R$.
\[ H(x) = \begin{cases}
            0, & \text{if $x < 0$,} \\
            \tfrac{1}{2}, & \text{if $x = 0$,} \\
            1, & \text{if $x > 0$.}
          \end{cases}
\]

\begin{lemma}[\cite{oesterle}*{Eq.~(37)}] \label{lem:des_cordes}
  Let $p$ be a prime and $k, k' \in \set{1, \ldots, p-1}$.
  We write $k_*$ for the unique element of~$\set{1, \ldots, p-1}$
  such that $k k_* \equiv -1 \bmod p$. Then
  \[ \xi(k) \bullet \xi(k') = -H(k'-k) + H(k'-k_*) + H(k'_*-k) - H(k'_*-k_*) . \]
\end{lemma}

\begin{proof}
  By~\cite{merel}*{Lemme~4}, for $k' \notin \set{k, k_*}$,
  $\xi(k) \bullet \xi(k')$ is the intersection number ($-1$, $0$, or~$1$)
  of the oriented line segment joining $e^{2\pi i k'_*/p}$ to~$e^{2\pi ik'/p}$
  and that joining $e^{2\pi i k_*/p}$ to~$e^{2\pi ik/p}$. Otherwise
  the intersection number is zero, since the pairing is alternating
  and $\xi(k_*) = -\xi(k)$. The formula we have given follows by considering
  the various possible cyclic orderings of the four points on the unit
  circle connected by the two line segments.
\end{proof}

We enlarge $M_n$ slightly and set
\[ B_n = \{(a,b,c,d) \in \Z^4 : a > b \ge 0, \; d > c \ge 0, \; ad - bc = n\} \]
(so we allow $c = 0$ here) and write $B_n^{b=0}$, $B_n^{b>0}$, $B_n^{c=0}$
and $B_n^{c>0} = M_n$ for the subsets satisfying the indicated extra
condition.

We define, for $n \ge 1$, a prime~$p > n$ and $k \in \set{1, \ldots, p-1}$,
the following two quantities.
\begin{align*}
  v_{p,n}(k) &= \#\{(a,b,c,d) \in \Z_{>0} : ad + bc = n, \; c \equiv dk \bmod p\} \\
  v'_{p,n}(k) &= \#\{(a,b,c,d) \in \Z_{>0} : ad + bc = n, \; \gcd(c,d) = 1, \; c \equiv dk \bmod p\}
\end{align*}

We now give an explicit formula for the intersection
number $(T_n - \sigma_1(n)) \bfe \bullet \xi(k)$. Its proof by Oesterl\'e
is quite ingenious.

\begin{proposition}[\cite{oesterle}*{Prop.~12}] \label{prop:In_inters_1}
  Let $p$ be a prime and $k, n \in \set{1, \ldots, p-1}$. Then
  \begin{enumerate}[\upshape(i)]
    \item \label{prop:In_inters_1:1}
          $\displaystyle (T_n - \sigma_1(n)) \bfe \bullet \xi(k)
             = \sum_{m \mid n} \Bigl(\floor{\frac{mk}{p}} - \floor{\frac{mk_*}{p}}\Bigr)
                 + v_{p,n}(k) - v_{p,n}(k_*)$.
    \item \label{prop:In_inters_1:2}
          $\displaystyle (T'_n - n) \bfe \bullet \xi(k)
             = \floor{\frac{nk}{p}} - \floor{\frac{nk_*}{p}}
                 + v'_{p,n}(k) - v'_{p,n}(k_*)$.
  \end{enumerate}
\end{proposition}

\begin{proof}
  Claim~\eqref{prop:In_inters_1:2} follows from
  claim~\eqref{prop:In_inters_1:1} by M\"obius inversion.
  So it suffices to show~\eqref{prop:In_inters_1:1}.

  We write $k_{c/d}$ for the integer $k \in \set{1, \ldots, p-1}$
  such that $c \equiv dk \bmod p$, where $c$ and~$d$ are integers
  coprime to~$p$. We extend this to all remaining elements $x \in \P^1(\Q)$
  by setting $k_x = p$. Then \cref{lem:Te_modsym,lem:des_cordes} imply that
  \begin{align*}
    (T_n &- \sigma_1(n)) \bfe \bullet \xi(k) \\
      &= -\sum_{(a,b,c,d) \in M_n} \xi(\tfrac{c}{d}) \bullet \xi(k) \\
      &= \sum_{(a,b,c,d) \in M_n}
           \bigl(H(k-k_{c/d}) - H(k-k_{-d/c}) - H(k_*-k_{c/d}) + H(k_*-k_{-d/c})\bigr) .
  \end{align*}
  We note that when $c=0$, all terms under the summation
  sign are zero (since $k_0 = k_\infty = p$ and $H(k - p) = H(k_* - p) = 0$),
  so that we can replace the summation over $M_n = B_n^{c>0}$ by a summation
  over~$B_n$ without changing the value of the sum.

  We now observe that there is a bijection
  \[ \phi_n \colon B_n^{b>0} \to B_n^{c>0}, \qquad
                   (a, b, c, d) \mapsto (b, -a+mb, d, -c+mb) ,
  \]
  where $m = \ceil{a/b} \ge 2$ is the unique integer such that
  $0 \le -a+mb < b$. (Its inverse is given by
  \[ (a, b, c, d) \mapsto (-b + m'a, a, -d + m'c, c) \qquad
     \text{with $m' = \ceil{d/c}$.)}
  \]
  We split the sum as follows.
  \begin{align}
    \sum_{(a,b,c,d) \in B_n} &\bigl(H(k - k_{c/d}) - H(k - k_{-d/c}))
                                      - H(k_*-k_{c/d}) + H(k_*-k_{-d/c})\bigr) \nonumber \\
      &= \sum_{(a,b,c,d) \in B_n^{b=0}} \bigl(H(k - k_{c/d}) - H(k_* - k_{c/d})\bigr)
         \nonumber \\
      &\qquad{} + \sum_{(a,b,c,d) \in B_n^{b>0}} \bigl(H(k - k_{c/d}) - H(k_* - k_{c/d})\bigr)
          \nonumber \\
      &\qquad{} - \sum_{(a,b,c,d) \in B_n^{c=0}} \bigl(H(k - k_{-d/c}) - H(k_* - k_{-d/c})\bigr)
          \nonumber \\
      &\qquad{} - \sum_{(a,b,c,d) \in B_n^{c>0}} \bigl(H(k - k_{-d/c}) - H(k_* - k_{-d/c})\bigr)
          \label{eq:imphi} .
  \end{align}
  Writing the quadruple in the sum~\eqref{eq:imphi} as~$\phi_n(a,b,c,d)$,
  this then gives
  \begin{align}
      &\sum_{(a,b,c,d) \in B_n^{b=0}} \bigl(H(k - k_{c/d}) - H(k_* - k_{c/d})\bigr)
         \label{eq:S1} \\
      &\quad{} - \sum_{(a,b,c,d) \in B_n^{c=0}} \bigl(H(k - k_{-d/c}) - H(k_* - k_{-d/c})\bigr)
         \label{eq:S2}\\
      &\quad{} + \sum_{(a,b,c,d) \in B_n^{b>0}}
                       \bigl(H(k - k_{c/d}) - H(k - k_{c/d - \ceil{a/b}}) \label{eq:S3} \\
      &\hspace*{35mm}{}      - H(k_* - k_{c/d}) + H(k_* - k_{c/d - \ceil{a/b}})\bigr). \nonumber
  \end{align}
  We evaluate the three sums in the last expression
  separately. First note that the second sum~\eqref{eq:S2} is zero,
  since $k_{-d/c} = p$ for $c = 0$ and $H(k - p) = H(k_* - p) = 0$ for all
  relevant~$k$. We now look at the first sum~\eqref{eq:S1}, which is the
  following expression minus the same expression with $k$ replaced by~$k_*$.
  \begin{align*}
    \sum_{(a,b,c,d) \in B_n^{b=0}} H(k - k_{c/d})
      &= \sum_{d \mid n} \sum_{c=0}^{d-1} H(k - k_{c/d})
       = \sum_{d \mid n} \sum_{c=1}^{d-1} H(k - k_{c/d}) .
  \end{align*}
  We set
  \begin{equation} \label{eq:def_sk}
    s(k) = \#\{(c,d) \in \Z^2 : d \mid n, \; 0 < c < d, \; c \equiv dk \bmod p\} ;
  \end{equation}
  then the sum above is
  \[ \sum_{d \mid n} \#\{c \in \Z : 0 < c < d, \; k_{c/d} \le k\} - \tfrac{1}{2} s(k) . \]
  Now $d k_{c/d} = up + c$, where $1 \le u < d$ satisfies
  $up \equiv -c \bmod d$, and so $k_{c/d} = \ceil{up/d}$. The map that sends
  $u$ to~$c$ is a permutation of~$\set{1, \ldots, d-1}$, which implies that
  \[ \#\{c \in \Z : 0 < c < d, \; k_{c/d} \le k\}
       = \#\{u \in \Z : 0 < u < d, \; u p \le d k\}
       = \floor{\frac{d k}{p}} .
  \]
  This gives the expression
  \begin{equation} \label{eq:sum1}
    \sum_{m \mid n} \left(\floor{\frac{m k}{p}} - \floor{\frac{m k_*}{p}}\right)
            - \tfrac{1}{2} \bigl(s(k) - s(k_*)\bigr)
  \end{equation}
  for the sum in~\eqref{eq:S1}.

  Now we look at the third sum~\eqref{eq:S3}.
  Let $x = c/d$ for some $(a,b,c,d) \in B_n^{b>0}$;
  then $p > n \ge d > 0$, so $p \nmid d$. If $\calA$ is a statement, we set
  $[\calA] = 0$ if $\calA$ is false and $[\calA] = 1$ if $\calA$ is true.
  Then, by checking the various cases and using that $k_{x-1} = k_x - 1$
  when $k_x \neq 1$, we find that
  \[ H(k - k_x) - H(k - k_{x-1}) = [k_x = 1] - \tfrac{1}{2} [k = k_x] - \tfrac{1}{2} [k = k_{x-1}] . \]
  This implies that
  \[ H(k - k_x) - H(k - k_{x-1}) - H(k_* - k_x) + H(k_* - k_{x-1})
      = \tfrac{1}{2} [k_* \in \{k_x, k_{x-1}\}] - \tfrac{1}{2} [k \in \{k_x, k_{x-1}\}] .
  \]
  We obtain the following expression for~\eqref{eq:S3}.
  \begin{align}
    \tfrac{1}{2}& \sum_{(a,b,c,d) \in B_n^{b>0}} \sum_{j=0}^{\ceil{a/b}-1}
                     \bigl([k_* \in \{k_{c/d-j}, k_{c/d-j-1}\}]
                             - [k \in \{k_{c/d-j}, k_{c/d-j-1}\}]\bigr) \nonumber \\
      &= \tfrac{1}{2}\bigl(\#\{(a,b,c,d) \in B_n^{b>0} : k_{c/d} = k_*\}
                              - \#\{(a,b,c,d) \in B_n^{b>0} : k_{c/d} = k\}\bigr)
        \label{eq:sum3_1} \\
      &\quad\begin{array}{@{}l}
              {} + \tfrac{1}{2}\bigl(\#\{(a,b,c,d) \in B_n^{b>0} : k_{c/d-\ceil{a/b}} = k_*\} \\[3pt]
              \qquad\quad{} - \#\{(a,b,c,d) \in B_n^{b>0} : k_{c/d-\ceil{a/b}} = k\}\bigr)
             \end{array} \label{eq:sum3_2} \\
      &\quad{} + \#\{(a,b,c,d,j) \in U_n : k_* = k_{c/d-j}\}
               - \#\{(a,b,c,d,j) \in U_n : k = k_{c/d-j}\} ,
        \label{eq:sum3_3}
  \end{align}
  where we have set
  \[ U_n = \{(a,b,c,d,j) : (a,b,c,d) \in B_n^{b>0}, \; 1 \le j < \ceil{a/b}\} . \]
  Now we observe that there is a bijection
  \[ \psi_n \colon U_n \to \{(a,b,c,d) \in \Z_{>0}^4 : ad + bc = n\}, \quad
                   (a,b,c,d,j) \mapsto (b, a-jb, d, -c+jd)
  \]
  (its inverse maps $(a,b,c,d)$ to $(b+ja, a, -d+jc, c, j)$ with
  $j = \ceil{d/c}$). Writing $\psi_n(a,b,c,d,j) = (a',b',c',d')$, we see
  that $k = k_{c/d-j}$ is equivalent to $k = k_{-d'/c'}$, which is the same
  as saying that $k_* = k_{c'/d'}$, or that $c' \equiv k_* d' \bmod p$.
  This shows that the terms in line~\eqref{eq:sum3_3} above are equal to
  \[ v_{p,n}(k) - v_{p,n}(k_*) . \]
  Using the bijection~$\phi_n$ between $B_n^{b>0}$ and~$B_n^{c>0}$,
  we see that the terms in line~\eqref{eq:sum3_2} can be written as
  \begin{align}
    \tfrac{1}{2}\bigl(\#\{(a,b,c,d) &\in B_n^{c>0} : k_{-d/c} = k_*\}
                          - \#\{(a,b,c,d) \in B_n^{c>0} : k_{-d/c} = k\}\bigr)
        \nonumber \\
    &= \tfrac{1}{2}\bigl(\#\{(a,b,c,d) \in B_n^{c>0} : k_{c/d} = k\}
                          - \#\{(a,b,c,d) \in B_n^{c>0} : k_{c/d} = k_*\}\bigr) .
        \label{eq:sum2_var}
  \end{align}
  This cancels the part of the terms in line~\eqref{eq:sum3_1} in
  which $c$ is strictly positive, and the terms with $c = 0$ do not
  contribute anything. What remains is the part with $b = 0$
  in~\eqref{eq:sum2_var}, which is
  \begin{align*}
    \tfrac{1}{2}\bigl(\#\{(c,d) &\in \Z_{>0}^2 : d > c > 0, \; d \mid n, \; c \equiv dk \bmod p\} \\
      &\quad{} - \#\{(c,d) \in \Z_{>0}^2 : d > c > 0, \; d \mid n, \; c \equiv dk_* \bmod p\}\bigr) \\
      &= \tfrac{1}{2}\bigl(s(k) - s(k_*)\bigr)
  \end{align*}
  with $s(k)$ as in~\eqref{eq:def_sk}. This cancels the contribution
  coming from $s(k)$ and~$s(k_*)$ in~\eqref{eq:sum1}, and we obtain the
  desired result.
\end{proof}

\begin{corollary} \label{cor:inters_large_p}
  Let $n \ge 1$ be an integer, let $c, d$ be coprime integers such
  that $c > d > 0$, and let $p > n c$ be a prime. Let $a, b$
  be the integers satisfying $0 \le a < c$, $0 \le b < d$ and $ad - bc = 1$.
  Let $k, k_* \in \{1,\ldots,p-1\}$ be such that $c \equiv d k \bmod p$
  and $-d \equiv c k_* \bmod p$. Further, let the integers $u$ and~$u_*$
  satisfy $dk = up + c$ and $ck_* = u_* p - d$.

  Then $0 \le u < d$, $0 \le u_* < c$, and
  \[ (T'_n - n) \bfe \bullet \xi(k)
       = \Bfloor{\frac{nu}{d}} - \Bfloor{\frac{nb}{d}}
          + \Bfloor{\frac{na}{c}} - \Bfloor{\frac{nu_*}{c}} .
  \]
\end{corollary}

\begin{proof}
  Since $dk - c > -p$ and $dk - c < dp$, we see that $0 \le u < d$.
  Since $ck_* + d > 0$ and $ck_* + d < c(k_* + 1) \le c p$, we also
  see that $0 \le u_* < c$.

  By~\cref{prop:In_inters_1},
  \[ (T'_n - n) \bfe \bullet \xi(k)
       = \Bfloor{\frac{nk}{p}} - \Bfloor{\frac{nk_*}{p}} + v'_{p,n}(k) - v'_{p,n}(k_*) .
  \]
  We evaluate each of the terms.

  We have that $nk/p = nu/d + nc/(pd)$ and $p > nc$, so $0 < nc/(pd) < 1/d$,
  which implies that $\floor{nk/p} = \floor{nu/d}$.

  Similarly, we have that $nk_*/p = nu_*/c - nd/(cp)$ and $p > nd$,
  so $0 < nd/(pc) < 1/c$, which implies that
  $\floor{nk_*/p} = \floor{(nu_* - 1)/c}$.

  The third term counts the quadruples $(a',b',c',d')$ of positive integers
  such that $c'$ and~$d'$ are coprime, $a'd' + b'c' = n$, and
  $c' \equiv d' k \bmod p$.
  The latter implies that $c' d \equiv c d' \bmod p$.
  Since $0 < c' d < n d < p$ and $0 < c d' < c n < p$, we must have
  equality; then the coprimality of $c'$ and~$d'$ and of $c$ and~$d$
  forces $(c',d') = (c,d)$.
  We have that $nad - nbc = n = a'd' + b'c' = a'd + b'c$, which implies
  that there is some $t \in \Z$ such that $na - a' = tc$ and $nb + b' = td$.
  The conditions $a', b' > 0$ then translate into $t < na/c$ and $t > nb/d$.
  Since $a/c > b/d$, this gives
  \begin{align*}
    v'_{p,n}(k)
      &= \#\Bigl\{t \in \Z : \frac{nb}{d} < t < \frac{na}{c}\Bigr\}
       = \floor{\frac{na-1}{c}} - \floor{\frac{nb}{d}} .
  \end{align*}
  The fourth term similarly counts quadruples $(a',b',c',d')$ of positive
  integers such that $c'$ and~$d'$ are coprime, $a'd' + b'c' = n$,
  and $p \mid c' k + d'$. The latter implies that $p \mid cc' + dd'$.
  But $0 < c c' + d d' < c (c' + d') \le c n < p$, so there are no such quadruples,
  and the fourth term is zero.

  Finally, note that
  \[ \Bfloor{\frac{na - 1}{c}} - \Bfloor{\frac{nu_* - 1}{c}}
      = \Bfloor{\frac{na}{c}} - \Bfloor{\frac{nu_*}{c}} ,
  \]
  as can be seen by considering the cases $c \mid n$ and $c \nmid n$
  separately, taking into account that $c$ is coprime with $a$ and~$u_*$.
\end{proof}

\begin{corollary} \label{cor:inters_aM}
  Let $M \ge 2$ be an integer, let $1 \le a < M$ be coprime with~$M$,
  let $n \ge 1$ be an integer, and let $p > n M$ be a prime.
  We let $w$ denote the integer such that $1 \le w < M$
  and $apw \equiv 1 \bmod M$. Then
  \[ (T'_n - n) \bfe \bullet \set{0, \frac{a}{M}}
      = \Bfloor{\frac{na}{M}} - \Bfloor{\frac{nw}{M}} .
  \]
\end{corollary}

\begin{proof}
  We prove this by induction on~$M$. When $M = 2$, then $a = 1$.
  We show the claim more generally for $a = 1$ and $M \ge 2$ arbitrary.
  We then have $\set{0, a/M} = \xi(M)$. The claim follows by taking
  $(c,d) = (M,1)$ (then $(a,b) = (1,0)$ and $(u,u_*) = (0,w)$)
  in~\cref{cor:inters_large_p}.

  Now assume that $M > 2$ and that the claim holds for smaller~$M$.
  We can then find integers $b$ and~$d$ such that $ad - bM = 1$
  and $1 \le d < M$. Then $0 \le b < d$. If $d = 1$, then $b = 0$
  and therefore $a = 1$; this case was already dealt with above.
  So we can assume that $d \ge 2$.

  Let $1 \le k < p$ be such that $M \equiv dk \bmod p$. Then
  \[ \begin{pmatrix} a-bk & b \\ M-dk & d \end{pmatrix}
       \cdot \set{0, \frac{1}{k}} = \set{\frac{b}{d}, \frac{a}{M}} .
  \]
  The matrix is in~$\Gamma_0(p)$, hence $\set{b/d, a/M} = \xi(k)$, so
  \[ (T'_n - n) \bfe \bullet \set{0, \frac{a}{M}}
       = (T'_n - n) \bfe \bullet \set{0, \frac{b}{d}}
          + (T'_n - n) \bfe \bullet \xi(k) .
  \]
  We use the induction hypothesis for the first term in the sum
  and~\cref{cor:inters_large_p} for the second term, where we take
  $(a,b,c,d) \leftarrow (a,b,M,d)$. Then
  \[ bpu = bdk - bc \equiv 1 \bmod d , \]
  so $u$ corresponds to~$w$
  in the induction hypothesis, and $u_* = w$. This gives
  \begin{align*}
    (T'_n - n) \bfe \bullet \set{0, \frac{a}{M}}
      &= \Bigl(\Bfloor{\frac{nb}{d}} - \Bfloor{\frac{nu}{d}}\Bigr)
            + \Bigl(\Bfloor{\frac{nu}{d}} - \Bfloor{\frac{nb}{d}}
                     + \Bfloor{\frac{na}{M}} - \Bfloor{\frac{nw}{M}}\Bigr) \\
      &= \Bfloor{\frac{na}{M}} - \Bfloor{\frac{nw}{M}} . \qedhere
  \end{align*}
\end{proof}

\begin{proof}[Proof of~\cref{lem:inters}]
  Using that $L_n = T'_{2n} - 2 T'_n = (T'_{2n} - 2n) - 2(T'_n - n)$,
  \cref{cor:inters_aM} gives (note that $w$ does not depend on~$n$)
  \begin{align*}
    L_n \bfe \bullet \set{0, \frac{a}{M}}
      &= \Bfloor{\frac{2na}{M}} - \Bfloor{\frac{2nw}{M}}
           - 2 \Bfloor{\frac{na}{M}} + 2 \Bfloor{\frac{nw}{M}}
       = \eps_M(na) - \eps_M(nw) ,
  \end{align*}
  and we can replace $w$ with $u/a$, where $u$ is as in~\cref{lem:inters}.
\end{proof}


\section{A criterion for ruling out moderately large primes} \label{sec:modlarge}

To exclude some of the larger primes for $d = 7$, we make use of the
following criterion, which is due to the first author of this paper.

\begin{proposition}[Derickx] \label{prop:hecke_gon}
  Let $d \ge 1$ and let $p$ be a prime. We assume that either
  \begin{enumerate}[\upshape(i)]
    \item \label{prop:hecke_gon_i}
          $J_1(p)(\Q)$ is finite, or
    \item \label{prop:hecke_gon_ii}
          there is $a \in (\Z/p\Z)^\times/\{\pm1\}$ such that $\ord(a) > 3d$
          and $A = (\diamondop{a}-1)(J_1(p)(\Q))$ is finite.
  \end{enumerate}
  In case~\eqref{prop:hecke_gon_ii},
  we say that ``$(*)$ holds'' when $\#A$ is odd or, more generally,
  the $2$-primary part of~$A$ is contained in the
  subgroup of~$J_1(p)(\Q)$ generated by differences of rational cusps.
  We then set $n = 3$ in case~\eqref{prop:hecke_gon_i} and
  \[ n = \begin{cases}
           5 & \text{if $(*)$ holds and $a \in \set{2, 2^{-1}}$,} \\
           6 & \text{if $(*)$ holds and $a \notin \set{2, 2^{-1}}$,} \\
           7 & \text{if $(*)$ does not hold and $a \in \set{3, 3^{-1}}$,} \\
           8 & \text{if $(*)$ does not hold and $a \notin \set{3, 3^{-1}}$.}
         \end{cases}
  \]
  in case~\eqref{prop:hecke_gon_ii}.
  Then $nd < \gon_{\Q}(X_1(p))$ implies that $p \notin S(d)$.
  This holds in particular when
  \[ d < \frac{325}{2^{16}} \, \frac{p^2 - 1}{n} \,. \]
\end{proposition}

\begin{proof}
  If $c \in X_1(p)$ is a rational cusp, which we consider as an effective
  divisor of degree~$1$, and $q$ is any prime, then
  $(T_q - \diamondop{q} - q)(c) = 0$.
  This can be deduced from the modular interpretation of the cusps.
  (See also~\cite{parent2}*{end of Section~2.4} and note that the rational
  cusps are those mapping to the cusp~$\infty$ on~$X_0(p)$.)

  We first consider case~\eqref{prop:hecke_gon_i}. Then
  by~\cref{cor:gen_by_rat_cusps}, $J_1(p)(\Q)$ is generated by differences
  of rational cusps. By the preceding paragraph, $T_q - \diamondop{q} - q$
  kills~$J_1(p)(\Q)$ for all primes~$q$ (including $q = 2$;
  this improves~\cref{prop:ann_rat_tors} in this case).
  In case~\eqref{prop:hecke_gon_ii}, \hbox{$T_q - \diamondop{q} - q$} kills
  the $2$-primary part of~$A$ when this is contained in the subgroup generated
  by differences of rational cusps and kills the odd part of~$A$
  by~\cref{prop:ann_rat_tors}. So when $(*)$ holds, $T_q - \diamondop{q} - q$
  kills~$A$ for arbitrary primes~$q$. When $(*)$ does not hold, the statement
  is true for $q \ge 3$.

  We let $x \in X_1(p)^{(d)}(\Q)$, considered as an effective divisor of
  degree~$d$ on~$X_1(p)$ and fix a rational cusp $c \in X_1(p)$.
  Then the linear equivalence class~$[x - d \cdot c]$ of the
  divisor $x - d \cdot c$ is a rational point on~$J_1(p)$.

  Going back to the case~\eqref{prop:hecke_gon_i},
  set $t = T_2 - \diamondop{2} - 2$. Then
  \[ t(x - d \cdot c) = t(x) - d t(c) = t(x) \]
  is a principal divisor, since $t([x - d \cdot c]) = 0$. This implies that
  the divisors $T_2(x)$ and $\diamondop{2}(x) + 2 x$ of degree $3d = nd$
  are linearly equivalent. But $\gon_{\Q}(X_1(p)) > nd$ by assumption, so
  the divisors must in fact be equal, and $t(x) = 0$. Now
  \cref{prop:hecke_op_kernel} shows that $x$ is a sum of cusps.
  This implies that $p \notin S(d)$.

  In case~\eqref{prop:hecke_gon_ii},
  we set $q = 2$ when $(*)$ holds and otherwise $q = 3$, so that
  $T_q - \diamondop{q} - q$ kills~$A$. Then $t(J_1(p)(\Q)) = \{0\}$, where
  \begin{equation} \label{eq:t_diff1}
    t = (\diamondop{a} - 1)(T_q - \diamondop{q} - q)
      = (\diamondop{a} T_q + \diamondop{q} + q)
          - (T_q + \diamondop{q a} + q \diamondop{a}) \,.
  \end{equation}
  If $qa = 1$, this simplifies to
  \begin{equation} \label{eq:t_diff2}
    t = (\diamondop{a} T_q + \diamondop{q} + (q-1)) - (T_q + q \diamondop{a}) \,,
  \end{equation}
  and if $a = q$, we obtain
  \begin{equation} \label{eq:t_diff3}
    t = (\diamondop{a} T_q + q) - (T_q + \diamondop{qa} + (q-1) \diamondop{a}) \,.
  \end{equation}
  We write $t_1$ for the first and~$t_2$ for the second term in the
  difference \eqref{eq:t_diff1}, \eqref{eq:t_diff2} or~\eqref{eq:t_diff3}.
  Since the diamond operators are automorphisms of~$X_1(p)$
  and $T_q$ multiplies degrees by~$q+1$, we see that
  applying~$t_1$ or~$t_2$, considered as a correspondence
  on~$X_1(p)$, to an effective divisor of degree~$d$ results in an
  effective divisor of degree~$nd$.

  As before, $t(x - d \cdot c) = t(x) - d t(c) = t(x)$ is a principal
  divisor, and we conclude from $\gon_{\Q}(X_1(p)) > nd$ that
  \[ t(x) = (T_q - \diamondop{q} - q)(\diamondop{a} - 1)(x) = 0 . \]
  By~\cref{prop:hecke_op_kernel} again, this implies that
  $\diamondop{a}(x) - x$ is supported on cusps. Since the diamond
  operators permute the cusps among themselves,
  this then implies that $x = x_0 + x_1$, where $x_0$ is supported in cusps
  and $x_1$ does not have cusps in its support and
  satisfies $\diamondop{a}(x_1) = x_1$. Now the diamond operators act
  freely on the non-cuspidal points of~$X_1(p)$ with the exception of
  points corresponding to elliptic curves with $j$-invariant $0$
  or~$1728$, which can have stabilizers of orders $3$ and~$2$, respectively.
  The condition $\diamondop{a}(x_1) = x_1$ implies that $x_1$ is a
  sum of (sums over) orbits of~$\diamondop{a}$, which have length
  at least $\ord(a)/3$. Since $\ord(a) > 3d$ by assumption, this forces
  $x_1 = 0$, and we conclude that $x$ is supported in cusps.
  This again implies that $p \notin S(d)$.

  For the last statement, note that
  \[ \gon_{\Q}(X_1(p)) \ge \gon_{\C}(X_1(p)) \ge \frac{\lambda_1}{48}(p^2 - 1) \]
  by~\cite{abramovich} (using that $\Gamma_1(p)$ has index $(p^2-1)/2$
  in $\PSL(2,\Z)$), where $\lambda_1$ is the smallest positive eigenvalue
  of the Laplace operator on~$X_1(p)(\C)$, which satisfies
  $\lambda_1 \ge 975/4096$ by~\cite{kim-sarnak}.
\end{proof}

\begin{remark}
  Without the condition ``$\ord(a) > 3d$'' in the case that $J_1(p)(\Q)$
  has positive rank, the proof shows that any rational point on~$X_1(p)^{(d)}$
  whose support consists of non-cuspidal points must be a sum of orbits
  of~$\diamondop{a}$. This is impossible when $d$ cannot be written
  as a sum of the possible orbit lengths (which are $\ord(a)$, together
  with $\ord(a)/2$ when $\ord(a)$ is even and $\ord(a)/3$ when $\ord(a)$
  is divisible by~$3$). But even when $d$ can be written in this way,
  this gives strong restrictions. For example, when $\ord(a) = d$
  and $d$ is coprime to~$6$, then such a point must be obtained by
  pulling back a rational point on~$X_H$, where $H$ is generated by~$a$.

  We plan to explore this further in a follow-up paper.
\end{remark}

\begin{corollary} \label{cor:7_71_113_127}
  We have that $p \notin S(7)$ for $p \in \{71, 113, 127\}$.
\end{corollary}

\begin{proof}
  We check that for the two primes $p \in \{113, 127\}$, the positive-rank
  simple factors of~$J_1(p)$ already occur in~$J_0(p)$. We can therefore
  take any $a \in (\Z/p\Z)^\times/\{\pm1\}$; we use $a = 3$, which
  generates $(\Z/p\Z)^\times/\{\pm1\}$ in both cases. In particular,
  $\ord(a) = (p-1)/2 > 3 \cdot 7$. We then have $n = 7$
  in~\cref{prop:hecke_gon}. Since
  \[ \frac{325}{2^{16}} \, \frac{p^2-1}{7} > 9 , \]
  all assumptions in~\cref{prop:hecke_gon} are satisfied.

  To deal with $p = 71$, we recall that by~\cref{prop:rank0}, $71$ is a
  rank zero prime, so we can apply~\cref{prop:hecke_gon} with $n = 3$. Since
  \[ \frac{325}{2^{16}} \, \frac{71^2-1}{3} > 8 , \]
  the claim follows also in this case.
\end{proof}

\begin{remark}
  For $p = 73$, the best we can do is use $a = 2$ and $n = 5$
  in~\cref{prop:hecke_gon} (by~\cite{CES}*{Section~6.2}, the torsion
  subgroup of~$J_1(73)(\Q)$ is generated by differences of rational
  cusps). However, the gonality lower bound works only for $d \le 5$.
  We would need $\gon_{\Q}(X_1(73)) > 35$. From
  Table~1 in~\cite{derickx_hoeij}, it appears that this is very likely
  the case, but it is also very likely hard to prove.
  (Note that $\ord(a) = 9 \le 3d$, but the argument would still work;
  see the remark following~\cref{prop:hecke_gon}.)
\end{remark}


\section{Verification of assumption \eqref{rc_b} of~\cref{lem:crit_residue_classes}} \label{sec:ass_b}

We now discuss assumption~\eqref{rc_b} of~\cref{lem:crit_residue_classes}
for the remaining pairs of degrees~$d$ and primes~$p$.
Recall that the assumption is always satisfied (with $\ell = 2$) when
\hbox{$p > (2^{d/2} + 1)^2$}; see~\cref{lem:alpha_surj_simple}.
The following table tells us which primes we still have to consider
for each~$d$.
\[ \renewcommand{\arraystretch}{1.25}
  \begin{array}{|r|ccccc|} \hline
                          d &  3 &  4 &  5 &  6 &   7 \\\hline
    \floor{(2^{d/2} + 1)^2} & 14 & 25 & 44 & 81 & 151 \\\hline
  \end{array}
\]

In some cases, we can show that all points in~$X_1(p)^{(d)}(\F_2)$
are sums of images of rational cusps, even when $p$ is below this bound.
The result of~\cite{waterhouse}*{Thm.~4.1} tells us precisely
what the possible orders of~$E(\F_{2^d})$ are for elliptic curves~$E$
defined over~$\F_{2^d}$. Using this (or a brute-force enumeration of
all such curves up to isomorphism), we obtain the following extension
of~\cref{lem:alpha_surj_simple}.

\begin{lemma} \label{lem:alpha_surj_ext}
  The set $X_1(p)^{(d)}(\F_2)$ consists of sums of images of rational cusps
  for $d$ and~$p$ as in the following table.
  \begin{align*}
    d = 3 \colon & \quad p = 11 \quad\text{or}\quad p \ge 17 \\
    d = 4 \colon & \quad p \ge 19 \\
    d = 5 \colon & \quad p \ge 23 \quad\text{and}\quad p \notin \set{31, 41} \\
    d = 6 \colon & \quad p = 23 \quad\text{or}\quad (p \geq 43 \;\;\text{and}\;\; p \neq 73) \\
    d = 7 \colon & \quad p \in \set{47, 53} \quad\text{or}\quad
                         (p \ge 79 \;\;\text{and}\;\; p \notin \set{113, 127})
  \end{align*}
\end{lemma}

\begin{proof}
  According to~\cite{waterhouse}*{Thm.~4.1}, $\#E(\F_{2^d})$ can take
  all even values in the Hasse interval
  $\left[\ceil{(2^{d/2}-1)^2}, \floor{(2^{d/2}+1)^2}\right]$
  and in addition the values
  \begin{align*}
    2^d + m 2^{d/2} + 1, & \quad m \in \{-2,-1,0,1,2\}, & & \text{if $d$ is even;} \\
    2^d + m 2^{(d+1)/2} + 1, & \quad m \in \{-1,0,1\}, & & \text{if $d$ is odd.}
  \end{align*}
  This allows us to determine the set of primes~$p$ such that there are no
  non-cuspidal points of degree~$\le d$ on~$X_1(p)_{\F_2}$.
  The condition that there are no cusps of degree~$\le d$ that are not
  images of rational cusps excludes in addition $p = 31$ for $d \ge 5$
  and $p = 127$ for $d \ge 7$.
\end{proof}

We note that for the primes not in the list above for a given~$d$,
there are indeed points~$\bar{x}$ as in assumption~\eqref{rc_b}.
If we want to show that $p \notin S(d)$ for one of these primes,
we have to do some work to show that there are no rational points
in the corresponding residue classes. For $p \in \set{29,31,41}$
and $d \ge 5$, we already did this in~\cref{lem:29_31_41_alpha}.
Taking into account~\cref{cor:7_71_113_127}, this leaves the primes
$p \in \set{37, 43, 59, 61, 67}$ for $d = 7$ and $p = 73$ for $d = 6, 7$.

We can deal with $(d,p) \in \set{(6,73), (7,43)}$ in the following way.

\begin{lemma} \label{lem:exclude_extra_primes}
  Let $d \ge 1$ be an integer and let $p$ be a prime.
  Let $\bar{x} \in X_1(p)^{(d)}(\F_2)$ be a point that is not a sum
  of images of rational cusps.
  Let $H \subseteq (\Z/p\Z)^\times/\set{\pm 1}$ be a subgroup
  and denote the image of~$\bar{x}$ in~$X_H^{(d)}(\F_2)$ by~$\bar{x}_H$.
  Assume that the following conditions are satisfied.
  \begin{enumerate}[\upshape(1)]
    \item \label{eep_1}
          There is $t \colon J_{H,\Z_{(2)}} \to A_{\Z_{(2)}}$ such that
          $t(J_H(\Q))$ is finite of odd order and  $t \circ \iota$
          (with $\iota \colon X_H^{(d)} \to J_H$) is a formal immersion
          at~$\bar{x}_H$.
    \item \label{eep_2}
          There is a rational point $x_H \in X_H^{(d)}(\Q)$ such that
          $\red_2(x_H) = \bar{x}_H$.
  \end{enumerate}
  Let $x \in X_1(p)^{(d)}(\Q)$ be such that $\red_2(x) = \bar{x}$.
  Then $x$ maps to~$x_H$ under the canonical map $X_1(p)^{(d)} \to X_H^{(d)}$.
\end{lemma}

\begin{proof}
  Let $x'_H$ be the image of~$x$ in~$X_H^{(d)}(\Q)$;
  then $\red_2(x'_H) = \bar{x}_H = \red_2(x_H)$.
  Since $t(J_H(\Q))$ is finite of odd order, this implies that
  $t(\iota(x'_H) - \iota(x_H)) = 0$.
  Since $t \circ \iota$ is a formal immersion at~$\bar{x}_H$,
  it follows that $x'_H = x_H$.
\end{proof}

If in the situation of~\cref{lem:exclude_extra_primes}, $x_H$ does not
lift to a rational point on~$X_1(p)^{(d)}$, then it follows that no rational
point on~$X_1(p)^{(d)}$ can reduce mod~$2$ to~$\bar{x}$.
We have to carry this out for all $\bar{x}$ as in assumption~\eqref{rc_b}.
To do this,
we formulate a criterion that allows us to verify the formal immersion
condition in~\cref{lem:exclude_extra_primes} also for points whose
support does not consist of cusps.

\begin{lemma} \label{lem:form_imm_crit}
  Fix a prime~$\ell$ and an integer $d \ge 1$.
  Let $X$ be a curve over~$\Q$ with good reduction at~$\ell$, with
  Jacobian variety~$J$. Fix $b \in X(\Q)$ and use it to define
  embeddings $\iota \colon X \to J$ and $\iota_d \colon X^{(d)} \to J$.
  Let $A$ be another abelian variety (with good reduction at~$\ell$)
  such that there is a homomorphism
  $t \colon J \to A$. Let $L \subseteq H^0(X_{\F_\ell}, \Omega^1)$
  be the pullback of~$H^0(A_{\F_\ell}, \Omega^1)$ under $t \circ \iota$,
  and let
  $\varphi \colon X_{\F_\ell} \to \P\Tan_0(A_{\F_\ell}) \cong \P^{\dim A - 1}_{\F_\ell}$
  be the morphism determined by the linear system corresponding to~$L$.
  Let $\bar{x} \in X^{(d)}(\F_\ell)$ be a point that is the sum of
  $d$~distinct geometric points
  $\bar{x}_1, \ldots, \bar{x}_d \in X(\bar{\F}_\ell)$. Assume that
  \begin{enumerate}[\upshape(i)]
    \item \label{fi_cond1}
          the differentials in~$L$ do not vanish simultaneously at any
          point~$\bar{x}_j$, and that
    \item \label{fi_cond2}
          the points
          $\varphi(\bar{x}_1), \ldots, \varphi(\bar{x}_d)
              \in \P^{\dim A - 1}(\bar{\F}_\ell)$
          span a linear subspace of dimension~$d-1$.
  \end{enumerate}
  Then $t \circ \iota_d$ is a formal immersion at~$\bar{x}$.
\end{lemma}

\begin{proof}
  To show that $t \circ \iota_d$ is a formal immersion, it is sufficient
  to show that the induced map on tangent spaces
  $\Tan_{\bar{x}}(X^{(d)}_{\F_\ell}) \to \Tan_{t(\iota(\bar{x}))}(A_{\F_\ell})$
  is injective; see~\cite{parent1}*{Thm.~4.18}.
  We can equivalently consider this condition over~$\bar{\F}_\ell$.

  Since the regular $1$-forms on~$A$ are invariant under translation,
  we have a canonical identification of all tangent spaces
  $\Tan_{\bar{a}}(A_{\bar{\F}_\ell})$ with the tangent space at the origin,
  whose projectivization is the codomain of~$\varphi$.
  Since the differentials in~$L$ do not vanish simultaneously at~$\bar{x}_j$,
  the map~$\varphi$ sends a point $\bar{x}_j \in X(\bar{\F}_\ell)$ to
  the image in~$\P\Tan_0(A_{\bar{\F}_\ell})$ of the tangent
  space~$\Tan_{\bar{x}_j}(X_{\bar{\F}_\ell})$
  under $(t \circ \iota)_*$ followed by a suitable translation.

  Since the geometric points making up~$\bar{x}$ are distinct, we have
  a canonical isomorphism
  \[ \Tan_{\bar{x}}(X^{(d)}_{\bar{\F}_\ell})
      \cong \bigoplus_{j=1}^d \Tan_{\bar{x}_j}(X_{\bar{\F}_\ell}) .
  \]
  The image of $\Tan_{\bar{x}}(X^{(d)}_{\bar{\F}_\ell})$
  in~$\P\Tan_0(A_{\bar{\F}_\ell})$
  under $(t \circ \iota_d)$ followed by a suitable translation is then
  the linear span of the various images~$\varphi(\bar{x}_j)$;
  the map on tangent spaces is injective if and only if
  this span has the maximal possible dimension~$d-1$.
\end{proof}

We will apply this as follows. We use the $q$-expansions mod~$2$
of the cusp forms associated to~$X_H$ to determine equations for
the canonical model of~$X_{H,\F_2}$. We then project away from the subspace
where the forms in~$L$ vanish (in practice, we compute the image
of~$\varphi$ in a similar way and then set up the projection) and
check that none of the points~$\bar{x}_j$ lie in this subspace. This
verifies the non-vanishing condition~\eqref{fi_cond1}. We then check
condition~\eqref{fi_cond2}.

\begin{lemma} \label{lem:6_73}
  Let $x \in X_1(73)^{(6)}(\Q)$. Then $\red_2(x) \in X_1(73)^{(6)}(\F_2)$
  is a sum of images of rational cusps.
\end{lemma}

\begin{proof}
  There are, up to isomorphism, exactly two elliptic curves over~$\F_{2^6}$
  with a point of order~$73$. They have zero $j$-invariant (they must
  be supersingular according to~\cite{waterhouse}) and automorphism
  group~$\Z/6\Z$, so each of them gives rise to $(73-1)/6 = 12$
  \hbox{$\F_{2^6}$-points} on~$X_1(73)_{\F_2}$. These $24$~points split
  into four orbits of size six under the action of Frobenius
  (each orbit contains three points coming from each of the two curves),
  so we obtain exactly four non-cuspidal points
  in~$X_1(73)^{(6)}(\F_2)$. There are no cuspidal points that are not
  sums of images of rational cusps, since the other cusps on~$X_1(73)_{\F_2}$
  are minimally defined over~$\F_{2^9}$. So we just have to exclude
  these four non-cuspidal points.

  Let $H$ be the subgroup
  of~$(\Z/73\Z)^\times/\set{\pm 1}$ of index~$9$. The canonical
  map $X_1(73) \to X_H$ is of degree~$4$ and unramified
  at all $24$~points mentioned above. This implies that they have six
  distinct images on~$X_H$; one can check that these points form one
  Frobenius orbit, so we get one point $\bar{x}_H \in X_H^{(6)}(\F_2)$
  that we have to deal with. The Jacobian~$J_H$ splits into a copy
  of~$J_{H'}$, where $H \subseteq H'$ has index~$3$,
  and a simple $30$-dimensional abelian variety~$A$.
  One can check that $A$ is a factor of the
  winding quotient and that all isogenous (over~$\Q$) abelian varieties
  have torsion subgroup of odd order (by computing orders
  of $A(\F_q)$ for suitable primes~$q$ via the Hecke eigenvalues).
  We take $t = \diamondop{7} - 1$; this kills~$J_{H'}$ and projects
  $J_H$ into~$A$.
  Since the nonzero eigenvalues of~$t$ are invertible mod~$2$
  (they are of the form $\omega - 1$ with $\omega \in \mu_3$),
  we can work with the $q$-expansions mod~$2$ of a basis of the
  space of cusp forms associated to~$A$. We check, as described above,
  that $t \circ \iota_6$ is a formal immersion at~$\bar{x}_H$.
  (In practice, we check this for all Frobenius orbits of length~$6$
  in $X_H(\bar{\F}_2)$, since it is not so easy to determine which
  point is in the support of~$\bar{x}_H$.)

  Note that $X_H \to X_{H'} \to X_0(73)$ is the composition of two
  maps of degree~$3$, the second of which is \'etale (by Riemann-Hurwitz:
  $X_{H'}$ is of genus~$13$ and $X_0(73)$ has genus~$5$).
  Let $E_0$ be an elliptic
  curve over~$\Q$ with complex multiplication by cube roots of unity.
  Then $E_0$ has two Galois-conjugate cyclic subgroups of order~$73$,
  with each subgroup defined over~$K = \Q(\sqrt{-3})$ (note that $73$
  splits in~$K$), so $E_0$ gives rise to a pair of Galois-conjugate
  points $y_1, y_2 \in X_0(73)(K)$. The preimages of these two points
  on~$X_{H'}$ give six geometric points that are Galois conjugate;
  the map $X_H \to X_{H'}$ is totally ramified at each of them, so we
  find a Galois orbit of size~$6$ of points in~$X_H(\bar{\Q})$,
  giving rise to a rational point $x_H \in X_H^{(6)}(\Q)$.
  This point reduces mod~$2$ to~$\bar{x}_H$ (as one can show by
  writing down an explicit twist of~$E_{0,K}$ for a certain number
  field of degree~$24$ that has a $K$-rational point of order~$73$
  and checking that the $24$~geometric points corresponding to
  its Galois conjugates reduce to the $24$ non-cuspidal points
  in~$X_1(73)(\F_{2^6})$ mentioned above), but does not lift to
  a rational point on~$X_1(73)^{(6)}$, since there are no CM elliptic
  curves with a $73$-torsion point over number fields of degree~$< 24$;
  see~\cite{CCS}*{Table~1}.
  By~\cref{lem:exclude_extra_primes} and the discussion following it,
  this finishes the proof.
\end{proof}

\begin{lemma} \label{lem:7_43}
  Let $x \in X_1(43)^{(7)}(\Q)$. Then $\red_2(x) \in X_1(43)^{(7)}(\F_2)$
  is a sum of images of rational cusps.
\end{lemma}

\begin{proof}
  There is, up to isomorphism, exactly one elliptic curve over~$\F_{2^7}$
  with a point of order~$43$. It is supersingular; its automorphism
  group over~$\F_{2^7}$ has order~$2$, since $\F_{2^7}$ does not contain
  primitive cube roots of unity. It therefore gives rise to~$21$ non-cuspidal
  points in~$X_1(43)(\F_{2^7})$, making up three Galois orbits.
  The non-rational cusps are also defined over~$\F_{2^7}$. We obtain
  six points in total in~$X_1(43)^{(7)}(\F_2)$ that are not supported
  in rational cusps. Take $H$ to be the subgroup of index~$7$.
  Then the six points above map to two points in~$X_H^{(7)}(\F_2)$.
  For~$A$, we use the winding quotient of~$J_H$; one can show that
  each $\Q$-isogenous abelian variety has odd torsion order.
  We show as before that $t \circ \iota$ is a formal
  immersion at the two points in question.

  On the other hand, there is a
  point in~$X_H^{(7)}(\Q)$ that corresponds to the pull-back of the
  cusp~$0$ on~$X_0(43)$ (note that $X_H \to X_0(43)$ has degree~$7$).
  It does not lift to a rational point on~$X_1(43)^{(7)}$, since the
  non-rational cusps on~$X_1(43)$ are points of degree~$21$. This
  shows that there are no rational points on~$X_1(43)^{(7)}$ whose
  reduction is cuspidal, but that are not supported in rational cusps.

  Consider now the rational point on~$X_0(43)$ that corresponds to
  elliptic curves over~$\Q$ with CM by the order of discriminant~$-43$.
  Its pullback to~$X_H$ again provides us with a rational point
  on~$X_H^{(7)}$, whose reduction must be the other point we have to
  consider, since such curves have (potentially) good reduction at~$2$.
  Again, this point does not lift to a rational point on~$X_1(43)^{(7)}$,
  as can be verified by consulting~\cite{CCS}*{Table~1}.
  This shows that there are no rational points on~$X_1(43)^{(7)}$ whose
  reduction is non-cuspidal.
\end{proof}

We still have to show that $p \notin S(7)$ for
\[ p = 37,\; 59,\; 61,\; 67,\; 73 . \]
We use the following simple observation by the first author
of this paper, together with the
fact that it is actually possible to check this criterion by a computation.

\begin{lemma}[Derickx] \label{lem:filter}
  Let $d \ge 1$ be an integer and let $p > 2$ be a prime. Assume that
  $t \in \T$ has the property that $t(J_1(p)(\Q)) = \{0\}$, where we
  consider $t$ as an endomorphism of~$J_1(p)$.
  Let $\bar{x}_0, \bar{x} \in X_1(p)^{(d)}(\F_2)$ be such that $\bar{x}_0$
  is a sum of images of rational cusps. If the divisor
  $t(\bar{x} - \bar{x}_0)$ on~$X_1(p)_{\F_2}$ is not principal
  (where we now consider $t$ as a correspondence on~$X_1(p)_{\F_2}$),
  then there is no rational point on~$X_1(p)^{(d)}$ whose reduction
  mod~$2$ is~$\bar{x}$.
\end{lemma}

\begin{remark}
  This result remains valid with an odd positive integer~$N$ in place
  of~$p$. (We need $N$ to be odd so that $X_1(N)$ has good reduction
  mod~$2$.)
\end{remark}

\begin{proof}
  Let $x_0 \in X_1(p)^{(d)}(\Q)$ be the sum of rational cusps such that
  $\red_2(x_0) = \bar{x}_0$ and assume that there is some
  $x \in X_1(p)^{(d)}(\Q)$ such that $\red_2(x) = \bar{x}$. Then
  the divisor $x - x_0$ represents a point in~$J_1(p)(\Q)$;
  it follows that $t(x - x_0)$ represents zero and is therefore principal.
  Applying reduction mod~$2$ shows that $t(\bar{x} - \bar{x}_0)$
  must be principal as well.
\end{proof}

We can find a suitable Hecke operator~$t$ by multiplying an operator
that projects $J_1(p)$ into an abelian subvariety of Mordell-Weil rank zero
(this is equivalent to this operator factoring through the winding quotient)
with an operator
that kills rational torsion. For the computations, we will use a
model of~$X_1(p)$ that is derived directly from the usual modular
interpretation, i.e., non-cuspidal points on~$X_1(p)$ correspond to pairs
$(E, P)$, where $E$ is an elliptic curve and $P \in E$ is a point
of exact order~$p$. The effect of a Hecke operator $T_n$ with $p \nmid 2n$
as a correspondence on~$X_1(p)_{\F_2}$ in this interpretation is then given
by mapping $(E, P)$ to the sum of the pairs $(E', \phi(P))$, where
$\phi \colon E \to E'$ runs through the cyclic isogenies of degree~$n$.
This switch from the ``natural'' modular interpretation given
in~\cref{sec:mod_curves} has the effect that we have to conjugate
everything by the Atkin-Lehner involution. Concretely, this means
that instead of $T_q - \diamondop{q} - q$ as stated
in~\cref{prop:ann_rat_tors}, we have to use $T_q - q \diamondop{q} - 1$
with any odd prime~$q$ to kill the rational torsion.
We will work with $q = 3$.

For the projection part of~$t$, we will use an operator of the form
$\diamondop{a} - 1$, so we take
\[ t = (\diamondop{a} - 1) (T_3 - 3 \diamondop{3} - 1) . \]
(This is similar to the idea used in~\cref{prop:hecke_gon}.)
We use the modular interpretation of the points on~$X_1(p)$ to find
the image of the divisor $\bar{x} - \bar{x}_0$ under~$t$.
Sutherland has computed planar equations for~$X_1(N)$ for all $N = p$ in
the relevant range, together with explicit expressions relating the
coordinates in these equations to the parameters $b$ and~$c$ in the
Tate form
\[ E_{b,c} \colon y^2 + (1-c) x y - b y = x^3 - b x^2 \]
of the associated elliptic curve with point $(0,0)$ of order~$N$.
See~\cite{sutherland}; the equations are available at~\cite{sutherland_table}.

We find the action of a diamond operator~$\diamondop{a}$ on a point
on~$X_1(p)$ by multiplying the point $P = (0,0)$ on the associated
curve~$E_{b,c}$ by~$a$ and then bringing the pair $(E_{b,c}, a P)$
into Tate form $(E_{b',c'}, (0,0))$. To get the effect of the Hecke
operator~$T_3$, we use the description of~$T_n$ given above, i.e.,
we find the four elliptic curves that are
$3$-isogenous to~$E_{b,c}$ (they may be defined over an extension
of the base field we are considering) and find the points corresponding
to the isogenous curves together with the image of~$P$. The sum of
these four points is then the image of the original point (considered
as a divisor of degree~$1$) under~$T_3$.

\begin{lemma} \label{lem:7_p}
  Let $p \in \{59, 61, 67, 73\}$ and $x \in X_1(p)^{(7)}(\Q)$.
  Then $\red_2(x) \in X_1(p)^{(7)}(\F_2)$ is a sum of images of
  rational cusps.
\end{lemma}

\begin{proof}
  We determine a suitable~$a$ for each of the primes~$p$ such
  that $\diamondop{a} - 1$ projects $J_1(p)$ into an
  abelian subvariety of rank zero.
  For $p \in \{59, 67, 73\}$, the only simple components of~$J_1(p)$
  that have positive rank are also components of~$J_0(p)$, so we
  can take $a$ to be any element of~$(\Z/p\Z)^\times/\{\pm1\}$.
  For $p = 61$, there
  is a component of positive rank in~$J_H$ for the subgroup~$H$ of
  index~$6$ that does not occur in~$J_0(p)$, and all components of
  positive rank occur in~$J_H$, so we take $a = 3 \equiv 2^6 \bmod 61$,
  where $2$ is a primitive root mod~$61$. We note that
  $\diamondop{a} - 1$ maps~$x_0$ to a degree zero divisor
  representing a torsion point in~$J_1(p)(\Q)$, so we just have to compute
  $t(\bar{x})$ and check whether this divisor is principal,
  where $\bar{x}$ and~$x_0$ are as in~\cref{lem:filter}.

  We then find all the non-cuspidal places of degree at most~$7$
  on~$X_1(p)_{\F_2}$. For the computation, it is sufficient to consider
  one representative in each orbit under the diamond operators.
  For $p < 73$, we find no such places of degree~$\le 6$ and either
  one (for $p = 61, 67$) or two (for $p = 59$) orbits of places
  of degree~$7$. For $p = 73$, there are two orbits of places of
  degree~$6$ and one orbit of places of degree~$7$.

  For the representatives~$\bar{x}$ of orbits of places of degree~$7$
  (which we identify with effective divisors of degree~$7$),
  we compute the divisor $t(\bar{x})$ and verify that it is not principal.
  This can be done by computing the Riemann-Roch space associated to
  the divisor; a divisor of degree zero is principal if and only if
  its Riemann-Roch space is nontrivial. (Magma has a built-in function
  for testing whether a divisor is principal.)

  The places of degree~$6$ on~$X_1(73)_{\F_2}$ give rise to effective
  divisors of degree~$7$ by adding one of the images of the rational
  cusps (which are exactly the $\F_2$-points on~$X_1(73)$).
  Applying~$t$ to such a sum differs from the result of applying~$t$ to
  the degree~$6$ divisor coming from the place by a principal divisor,
  since the rational cusps map to principal divisors. So we only
  have to check that $t(\bar{x})$ is non-principal for the two
  representatives of orbits of places of degree~$6$.
  (We note that this also gives an alternative proof of~\cref{lem:6_73}.)

  Finally, we note that all other points in $X_1(p)^{(7)}(\F_2)$
  are supported in images of rational cusps, since the other cusps
  give rise to points of degree at least~$9$ over~$\F_2$.

  The computations took less than one hour each for $p = 59$ and~$61$,
  about three hours for $p = 67$ and about seven hours for $p = 73$.
\end{proof}

\begin{remark}
  We note that we can use this approach also to show that there are
  no non-cuspidal points in~$X_1(43)^{(7)}(\F_2)$ that arise as the
  reduction modulo~$2$ of a rational point. We would still have to
  deal with the points arising from Frobenius orbits of cusps that
  are not images of rational cusps, however; see the proof
  of~\cref{lem:7_43}.
\end{remark}

Now \cref{prop:ass_b} follows
from~\cref{lem:alpha_surj_ext,lem:29_31_41_alpha,cor:7_71_113_127,lem:6_73,lem:7_43,lem:7_p}.

Finally, we deal with $p = 37$.

\begin{lemma} \label{lem:37_6}
  Modulo the action of Frobenius and the diamond operators, there is
  exactly one point of degree~$6$ on~$X_1(37)_{\F_2}$ such that the
  corresponding point $\bar{x} \in X_1(37)^{(6)}(\F_2)$ is the reduction
  mod~$2$ of a rational point $x \in X_1(37)^{(6)}(\Q)$, and this point~$x$
  is uniquely determined by~$\bar{x}$.
\end{lemma}

\begin{proof}
  We proceed as in the proof of~\cref{lem:7_p}. The only positive-rank
  factor of~$J_1(37)$ occurs in~$J_0(37)$ (it is the ``first'' elliptic
  curve of rank~$1$), so we can take any~$a$ for the criterion
  of~\cref{lem:filter}.
  The computation shows that of the two diamond orbits of places of
  degree~$6$, only one satisfies the criterion in~\cref{lem:filter}.
  (It should be noted that this can be used to verify that we are
  correct in working with the Hecke operator $T_3 - 3 \diamondop{3} - 1$:
  none of the two places satisfies the criterion when using
  $T_3 - \diamondop{3} - 3$ instead, but one of them has to, since
  there are non-cuspidal rational points on~$X_1(37)^{(6)}$.)

  We know that there is a diamond orbit of
  rational points that has to reduce to our unique diamond orbit that
  lifts. To show that the lift is unique, we use~\cref{lem:form_imm_crit}.
  The Hecke operator~$T_{17}$ projects $J_1(37)$ into an abelian
  subvariety of rank zero. Its eigenvalues are invertible mod~$2$
  on newforms corresponding to a subvariety of dimension~$36$,
  which has odd order rational torsion subgroup. We then verify the formal
  immersion criterion (for all points of degree~$6$, since we work with
  a different model here and did not try to find an explicit birational
  map between the two models).
\end{proof}

\begin{proof}[Proof of~\cref{prop:6_37}]
  Let $x \in X_1(37)^{(6)}(\Q)$ be a point whose support contains
  no cusps. Since~\eqref{rc_a} holds for $(d, p) = (6, 37)$
  by~\cref{prop:ass_a} and there are no non-cuspidal points
  on~$X_1(37)_{\F_2}$ of degree~$\le 5$, it follows that
  $\bar{x} = \red_2(x) \in X_1(37)^{(6)}(\F_2)$ is also a point
  whose support contains no cusps. By~\cref{lem:37_6}, $\bar{x}$
  is uniquely determined up to the action of the diamond operators,
  and there is no other point than~$x$ that reduces mod~$2$ to~$\bar{x}$.
  On the other hand, we know a point~$x'$ with this property; this is a
  point coming from the curve~$E_{6,37}$ with some choice of point
  of order~$37$ (they are all in the same diamond orbit).
  It follows that $x = x'$, which implies the claim.
\end{proof}

We finish off the determination of~$S(7)$ by excluding $p = 37$.

\begin{lemma} \label{lem:37_7}
  $37 \notin S(7)$.
\end{lemma}

\begin{proof}
  As in the proof of~\cref{lem:37_6}, we show that there is no point
  of degree~$7$ on~$X_1(37)_{\F_2}$ such that the corresponding
  point in $X_1(37)^{(7)}(\F_2)$ is the reduction of a rational point.
  Now assume that $x \in X_1(37)^{(7)}(\Q)$ and consider
  $\bar{x} = \red_2(x)$. By the preceding statement, the support
  of~$\bar{x}$ must contain a cusp, and the non-cuspidal part
  of~$\bar{x}$ must satisfy the criterion of~\cref{lem:filter}.
  By~\cref{lem:37_6} and its proof, the non-cuspidal part is then
  either empty, or it is in the unique diamond orbit coming from
  non-cuspidal rational points on~$X_1(37)^{(6)}$. In the first
  case, $x$ must be a sum of rational cusps, since
  assumption~\eqref{rc_a} holds. To deal with the second case,
  we verify the formal immersion criterion as in the proof
  of~\cref{lem:37_6}, but now for all sums of an $\F_2$-rational cusp
  and a prime divisor of degree~$6$. This shows that the criterion
  is satisfied; therefore the point~$x$ is unique in its residue
  class mod~$2$. On the other hand, there is a known point in this
  residue class, which comes from adding the rational cusp that lifts
  the unique cusp in the support of~$\bar{x}$ to the degree~$6$ divisor
  lifting the remaining part (this is one of the sporadic points
  in~$X_1(37)^{(6)}(\Q)$). It follows that $x$ is this point; in
  particular, $x$ has a cusp in its support. So we conclude that every
  rational point on~$X_1(37)^{(7)}$ has a cusp in its support;
  this is equivalent to the statement that $37 \notin S(7)$.
\end{proof}


\end{document}